\documentclass[a4paper]{article}

\usepackage{arydshln}
\usepackage{amsmath}
\usepackage{mathtools}
\usepackage{amsfonts,amssymb}
\usepackage[compress]{cite}
\usepackage{graphicx}
\usepackage{mdwlist}
\usepackage{xfrac}
\usepackage{paralist}
\usepackage{algorithmic}
\usepackage{algorithm}
\usepackage{color}
\usepackage{url}

\newtheorem{lemma}{\sc Lemma}[section]
\newtheorem{theorem}[lemma]{\sc Theorem}

\newtheorem{assumption}[lemma]{\sc Assumption}

\newtheorem{corollary}[lemma]{\sc Corollary}

\renewcommand{\matrix}[2]{\left[\begin{array}{#1} #2 \end{array}\right] }

\DeclareMathOperator*{\trace}{trace}
\DeclareMathOperator*{\diag}{diag}

\usepackage[font=footnotesize]{caption}

\newcommand{\myQED}{~\rule[-1pt]{5pt}{5pt}\par\medskip}
\newenvironment{myproof}{{\noindent \bf Proof:\ }}{ \hfill \myQED}

\begin{document}

\title{Stochastic Sensor Scheduling for \\ Networked Control Systems\thanks{A preliminary and brief version of this work was presented at the American Control Conference, 2013~\cite{FarokhiSubmittedACC2012}. The work was supported by the Swedish Research Council and the Knut and Alice Wallenberg Foundation.}}
\author{Farhad~Farokhi$^\dag$~and~Karl~H.~Johansson\thanks{The authors are with ACCESS Linnaeus Center, School of Electrical Engineering, KTH Royal Institute of Technology, Stockholm, Sweden. E-mails: \{farokhi,kallej\}@ee.kth.se }}

\date{}

\maketitle

\begin{abstract} Optimal sensor scheduling with applications to networked estimation and control systems is considered. We model sensor measurement and transmission instances using jumps between states of a continuous-time Markov chain. We introduce a cost function for this Markov chain as the summation of terms depending on the average sampling frequencies of the subsystems and the effort needed for changing the parameters of the underlying Markov chain. By minimizing this cost function through extending Brockett's recent approach to optimal control of Markov chains, we extract an optimal scheduling policy to fairly allocate the network resources among the control loops. We study the statistical properties of this scheduling policy in order to compute upper bounds for the closed-loop performance of the networked system, where several decoupled scalar subsystems are connected to their corresponding estimator or controller through a shared communication medium. We generalize the estimation results to observable subsystems of arbitrary order. Finally, we illustrate the developed results numerically on a networked system composed of several decoupled water tanks. 
\end{abstract}

\begin{figure}[t]
\centering
$$
\begin{array}{p{.50\linewidth}p{.07\linewidth}p{.22\linewidth}}
\includegraphics[width=0.9\linewidth]{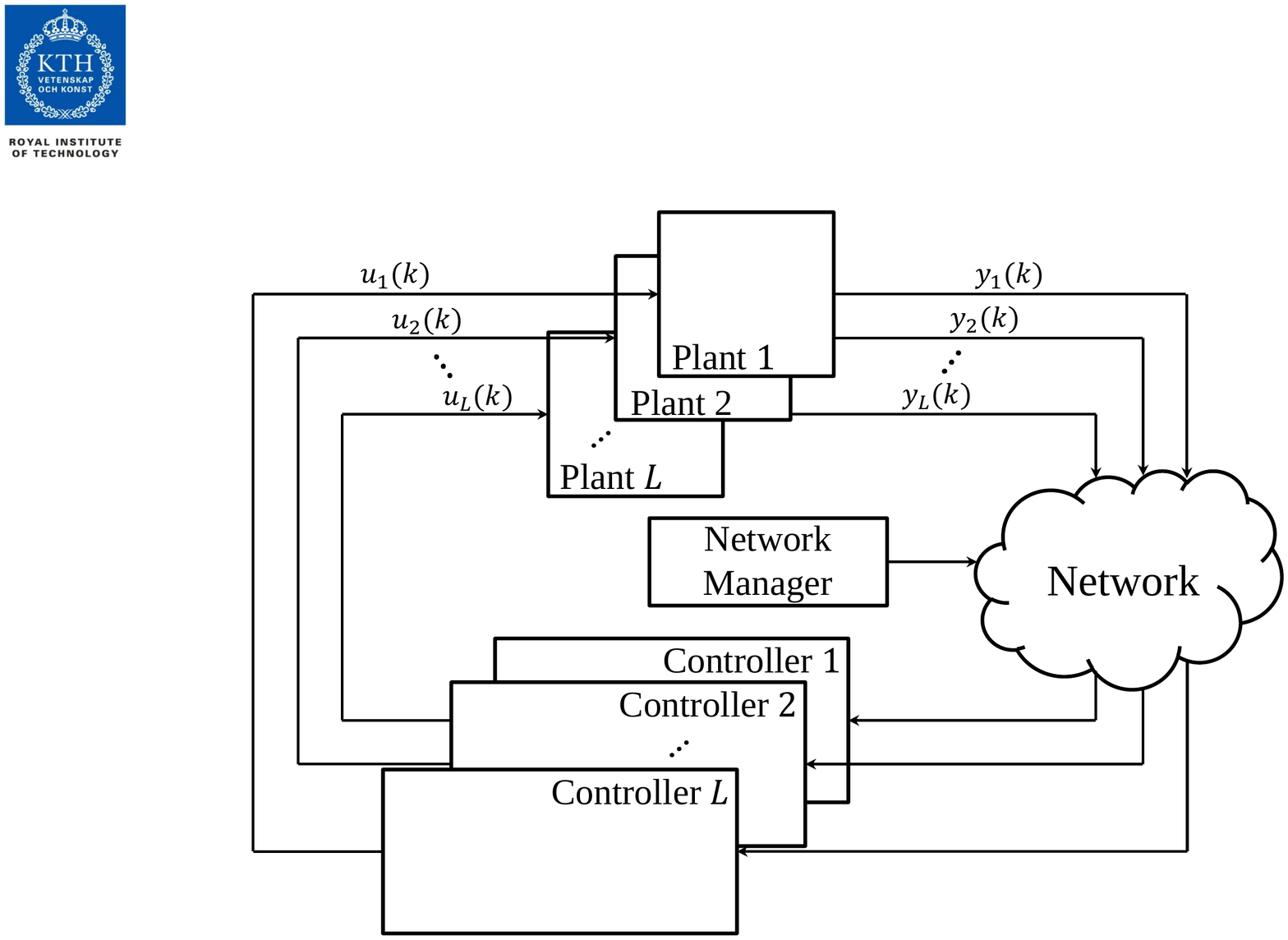} & \hspace{.07\linewidth} &\includegraphics[width=0.9\linewidth]{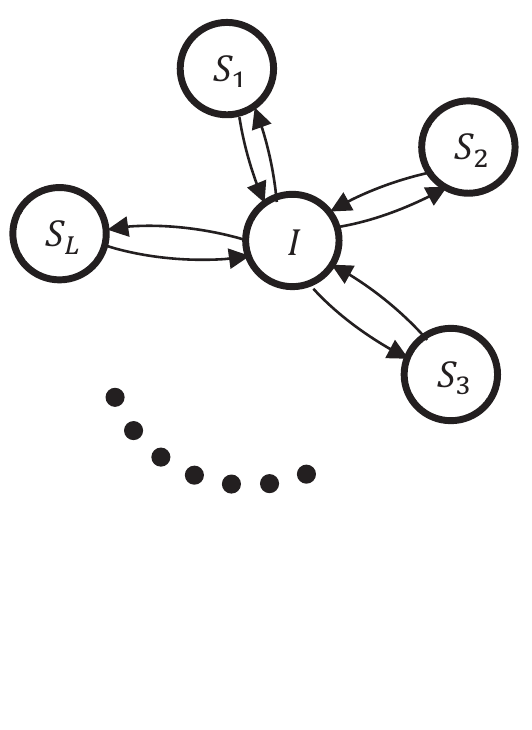}
\end{array}
$$
\caption{\label{figureNCS_new} An example of a networked control system (left) together with flow diagram of the continuous-time Markov chain used for modeling the proposed stochastic scheduling policy (right).}
\end{figure}

\section{Introduction}
\subsection{Motivation}
Emerging large-scale control applications in smart infrastructures~\cite{negenborn2010intelligent}, intelligent transportation systems~\cite{Varaiya250509}, aerospace systems~\cite{Giulietti887447}, and power grids~\cite{Massoud1507024}, are typically implemented over a shared communication medium. Figure~\ref{figureNCS_new} illustrates an example of such a networked system, where~$L$ decoupled subsystems are connected to their subcontrollers over a wireless communication network. A set of sensors in each subsystem sample its state and transmit the measurements over the wireless network to the corresponding subcontroller. Then, the subcontroller calculates an actuation signal (based on the transmitted observation history) and directly applies it to the subsystem. Unfortunately, traditional digital control theory mostly results in conservative networked controllers because the available methods often assume that the sampling is done periodically with a fixed rate~\cite{astrom1984computer,franklin1998digital}. When utilizing these periodic sampling methods, the network manager should allocate communication instances (according to the fixed sampling rates) to each control loop considering the worst-case possible scenario, that is, the maximum number of active control loops. In a large control system with thousands of control loops, fixed scheduling of communication instances imposes major constraints because network resources are allocated even if a particular control loop is not active at the moment. This restriction is more evident in ad-hoc networked control systems where many control loops may join or leave the network or switch between active and inactive states. Therefore, we need a scheduling method to set the sampling rates of the individual control loops adaptively according to their requirements and the overall network resources. We address this problem in this paper by introducing an optimal stochastic sensor scheduling scheme. 

\subsection{Related Studies}
In networked control systems, communication resources need to be efficiently shared between multiple control loops in order to guarantee a good closed-loop performance. Despite that communication resources in large networks almost always are varying over time due to the need from the individual users and physical communication constraints, the early networked control system literature focused on situations with fixed communication constraints; e.g., bit-rate constraints~\cite{Elia948466,tatikonda2004control,Nair2003585,Brockett763226} and packet loss~\cite{Schenato4118476,Sinopoli1333199,Smith1254085,gupta2007optimal}. Only recently, some studies have targeted the problem of integrated resource allocation and feedback control; e.g.,~\cite{ramesh2011dual,ramesh2009lqg,molin2009lqg,quagli2010design,nilsson1997lqg,fontanelli2013soft}.

The problem of sharing a common communication medium or processing unit between several users is a well-known problem in computer science, wireless communication, and networked control~\cite{Walsh898792,2006resource,buttazzo2005hard,liu2006firm}. For instance, the authors in~\cite{kay1988fair} proposed a scheduler to allocate time slots between several users over a long horizon. In that scheduler, the designer must first manually assign shares (of a communication medium or processing unit) that an individual user should receive. Then, each user achieves its pre-assigned share by means of probabilistic or deterministic algorithms~\cite{kay1988fair,waldspurger1994lottery}. The authors in~\cite{jackson1955scheduling,horn1974some} proved that implementing the task with the earliest deadline achieves the optimum latency in case of both synchronous and asynchronous job arrivals. In~\cite{liu1973scheduling}, a scheduling policy based on static priority assignment to the tasks was introduced. Many studies in communication literature have also considered the problem of developing protocols in order to avoid the interference between several information sources when using a common communication medium. Examples of such protocols are both time-division and frequency-division multiple access~\cite{goldsmith2005wireless,tse2005fundamentals}. Contrary to all these studies, in this paper, we automatically determine the communication instances (and, equivalently, the sampling rates) of the subsystems in a networked system based on the number of active control loops at any given moment. We use a continuous-time Markov chain to model the optimal scheduling policy.

Markov chains are very convenient tools in control and communication~\cite{shannon1962mathematical,meyn2007control}. Markov jump linear systems with underlying parameters switching according to a given Markov chain has been studied in the control literature~\cite{Lee2006205,feng1992stochastic,costa2004discrete,blair1975feedback}. The problem of controlled Markov chains has always been actively pursued~\cite{bellman2003dynamic,howard1960dynamic,watkins1992q,meyn2009markov}. In a recent study by Brockett~\cite{4738725}, an explicit solution to the problem of optimal control of observable continuous-time Markov chains for a class of quadratic cost functions was presented. In that paper, the underlying continuous-time Markov chain was described using the so-called unit vector representation~\cite{4738725,LectureA}. Then, the finite horizon problem and its generalization to infinite horizon cost functions were considered. We extend that result to derive the optimal scheduling policy in this paper.

In the study~\cite{gupta2006stochastic}, the authors developed a stochastic sensor scheduling policy using Markov chains. Contrary to this paper, they considered a discrete-time Markov chain to get a numerically tractable algorithm for optimal sensor scheduling. The algorithm in~\cite{gupta2006stochastic} uses one of the sensors at each time step while here, the continuous-time Markov chain can rest in one of its states to avoid sampling any of the sensors. Furthermore, the cost function in~\cite{gupta2006stochastic} was not written explicitly in terms of the Markov chain parameters, but instead it was based on the networked system performance when using a Markov chain for sampling the sensors. However, our proposed scheduling policy results in a separation between designing the Markov chain parameters and networked system, which enables us to describe the cost function needed for deriving our optimal sensor scheduling policy only in terms of the Markov chain parameters. 

\subsection{Main Contributions}
The objective of the paper is to find a dynamic scheduling policy to fairly allocate the network resources between the subsystems in a networked system such as the one in Figure~\ref{figureNCS_new}~(left). Specifically, we employ a continuous-time Markov chain for scheduling the sensor measurement and transmission instances. We use time instances of the jumps between states of this continuous-time Markov chain to model the sampling instances; i.e., whenever there is a jump from an idle state in the Markov chain to a state that represent a subsystem in the networked system, we sample that particular subsystem and transmit its state measurement across the shared communication network to the corresponding subcontroller. Figure~\ref{figureNCS_new}~(right) illustrates the flow diagram of the proposed Markov chain. Every time that a jump from the idle node $I$ to node $S_\ell$, $1\leq \ell \leq L$, occurs in this continuous-time Markov chain, we sample subsystem $\ell$ and send its state measurement to subcontroller~$\ell$. The idle state $I$ helps to tune the sampling rates of the subsystems independently. As an approximation of the wireless communication network, we assume that the sampling and communication are instantaneous; i.e., the sampling and transmission delays are negligible in comparison to the subsystems response time. We still want to limit the amount of communication per time unit to reduce the energy consumption and network resources.

We mathematically model the described continuous-time Markov chain using unit vector representation~\cite{4738725,LectureA}. We introduce a cost function that is a combination of the average sampling frequencies of the subsystems (i.e., the average frequency of the jumps between the idle state and the rest of the states in the Markov chain) and the effort needed for changing the scheduling policy (i.e., changing the underlying Markov chain parameters). We expand the results presented in~\cite{4738725} to minimize the cost function over both finite and infinite horizons. Doing so, we find an explicit minimizer of the cost function and develop the optimal scheduling policy accordingly. This policy fairly allocates sampling instances among the sensors in the networked system. The proposed optimal scheduling policy works particularly well for ad-hoc sensor networks since we can easily accommodate for the changes in the network configuration by adding an extra state to the Markov chain (and, in turn, by adding an extra term to the cost function) whenever a new sensor becomes active and by removing a state from the Markov chain (and, in turn, by removing the corresponding term from the cost function) whenever a sensor becomes inactive. The idea of dynamic peer participation (or churn) in peer-to-peer networks have been extensively studied in the communication literature~\cite{stutzbach2006understanding,shen2007locality}. However, not much attention has been paid to this problem for networked control and estimation.

Later, we focus on networked estimation as an application of the proposed stochastic sensor scheduling policy. We start by studying a networked system composed of several scalar subsystems and calculate an explicit upper bound for the estimation error variance as a function of the statistics of the measurement noise and the scheduling policy. The statistics of the scheduling policy are implicitly dependent on the cost function. Hence, we can achieve the required level of performance by finely tuning the cost function design parameters. We generalize these estimation results to higher-order subsystems when noisy state measurements of the subsystems are available. In the case where noisy output measurements of the subsystems are available, we derive an estimator based on the discrete-time Kalman filter and calculate an upper bound for the variance of its error given a specific sequence of sampling instances. Lastly, we consider networked control as an application of the proposed sensor scheduling policy. We assume that the networked control system is composed of scalar subsystems that are in feedback interconnection with impulsive controllers (i.e., controllers that ideally reset the state of the system whenever a new measurement arrives). We find an upper bound for the closed-loop performance of the subsystems as a function of the statistics of the measurement noise and the scheduling policy. We generalize this result to pulse and exponential controllers. 

\subsection{Paper Outline}
The rest of the paper is organized as follows. In Section~\ref{sec:NetConSys}, we introduce the optimal stochastic scheduling policy and calculate its statistics. We apply the proposed stochastic scheduling policy to networked estimation and control systems in Sections~\ref{sec:NetworkedEstimationPerformance} and~\ref{sec:NetworkedControlPerformance}, respectively. In Section~\ref{sec:NumericalExample}, we illustrate the developed results numerically on a networked system composed of several decoupled water tanks. Finally, we present the conclusions and directions for future research in Section~\ref{sec:Conclusion}. 

\subsection{Notation}
The sets of integer and real numbers are denoted by $\mathbb{Z}$ and $\mathbb{R}$, respectively. We use $\mathbb{O}$ and $\mathbb{E}$ to denote the sets of odd and even numbers. For any $n\in\mathbb{Z}$ and $x\in\mathbb{R}$, we define $\mathbb{Z}_{> (\geq) n}=\{m\in\mathbb{Z}\;|\; m> (\geq) n\}$ and $\mathbb{R}_{> (\geq) x}=\{y\in\mathbb{R}\;|\; y> (\geq) x\}$, respectively. We use calligraphic letters, such as $\mathcal{A}$ and $\mathcal{X}$, to denote any other set.

We use capital roman letters, such as $A$ and $C$, to denote matrices. For any matrix $A$, $a_{ij}$ denotes its entry in the $i$-th row and the $j$-th column.

Vector $e_i$ denotes a column vector (where its size will be defined in the text) with all entries equal zero except its $i$-th entry which is equal to one. For any vector $x\in\mathbb{R}^n$, we define the entry-wise operator $x^{.2}=[x_1^2\;\; \cdots\;\;x_n^2]^\top$.

\section{Stochastic Sensor Scheduling} \label{sec:NetConSys}
In this section, we develop an optimal stochastic scheduling policy for networked systems, where several sensors are connected to the corresponding controllers or estimators over a shared communication medium. Let us start by modeling the stochastic scheduling policy using con-tinuous-time Markov chains. 

\subsection{Sensor Scheduling Using Continuous-Time Markov Chains} \label{subsec:ssns}
We employ continuous-time Markov chains to model the sampling instances of the subsystems. To be specific, every time that a jump from the idle node $I$ to node $S_\ell$, $1\leq \ell \leq L$, occurs in the continuous-time Markov chain described by the schematic flow diagram in Figure~\ref{figureNCS_new}~(right), we sample subsystem $\ell$. We use unit vector representation to mathematically model this continuous-time Markov chain~\cite{4738725,LectureA}.

We define the set $\mathcal{X}=\{e_1,e_2,\dots,e_{n}\}\subset \mathbb{R}^n$ where $n=L+1$. The Markov chain state $x(t)\in\mathbb{R}^n$ takes value from $\mathcal{X}$, which is the reason behind naming this representation as the unit vector representation. We associate nodes $S_1$, $S_2$, $\dots$, $S_L$, and $I$ in the Markov chain flow diagram with unit vectors $e_1$, $e_2$, $\dots$, $e_{L}$, and $e_n$, respectively. Following the same approach as in~\cite{LectureA}, we can model the Markov chain in Figure~\ref{figureNCS_new}~(right) by the It\^{o} differential equation
\begin{equation} \label{eqn:uvr:1}
\mathrm{d}x(t)=\sum_{\ell=1}^L \bigg(G'_{\ell n} x(t)\;\mathrm{d}N'_{\ell n}(t)+G'_{n\ell}x(t)\;\mathrm{d}N'_{n\ell}(t)\bigg), 
\end{equation}
where $\{N'_{n\ell}(t)\}_{t\in\mathbb{R}_{\geq 0}}$ and $\{N'_{\ell n}(t)\}_{t\in\mathbb{R}_{\geq 0}}$, $1\leq \ell \leq L$, are Poisson counter processes\footnote{Recall that a Poisson counter $N(t)$ is a stochastic process with independent and stationary increments that starts from zero $N(0)=0$. Additionally, $\mathbb{P}\{N(t+\Delta t)-N(t)=k\}=(\int_{t}^{t+\Delta t} \lambda(t) \mathrm{d}t)^k\exp(\int_{t}^{t+\Delta t} \lambda(t) \mathrm{d}t)/k!$ for any $t,\Delta t\in\mathbb{R}_{\geq 0}$ and $k\in\mathbb{Z}_{\geq 0}$. In the limit, when replacing $\Delta t$ with $\mathrm{d}t$ and $\Delta N(t)=N(t+\Delta t)-N(t)$ with $\mathrm{d}N(t)$, we get $\mathbb{P}\{\mathrm{d}N(t)=0\}=1-\lambda(t)\mathrm{d}t$, $\mathbb{P}\{\mathrm{d} N(t)=1\}=\lambda(t)\mathrm{d}t$, and $\mathbb{P}\{\mathrm{d}N(t)=k\}=0$ for $k\in\mathbb{Z}_{\geq 2}$. For a detailed discussion on Poisson counters see~\cite{cox1980point,LectureA}.} with rates $\lambda_{n\ell}(t)$ and $\lambda_{\ell n}(t)$, respectively. These Poisson counters determine the rates of jump from $S_\ell$ to~$I$, and vice versa. In addition, we have $G'_{\ell n}=(e_\ell-e_n)e_n^\top$ and $G'_{n\ell }=(e_n-e_\ell)e_\ell^\top$, $1\leq \ell \leq L$. Let us define $m=2L$. Now, we can rearrange the It\^{o} differential equation in~(\ref{eqn:uvr:1}) as
\begin{equation} \label{eqn:uvr:sys}
\mathrm{d}x(t)=\sum_{i=1}^{m} G_i x(t)\;\mathrm{d}N_i(t),
\end{equation}
where $\{N_i(t)\}_{t\in\mathbb{R}_{\geq 0}}$, $1\leq i\leq m$, is a Poisson counter process with rate denoted as
\begin{equation} \label{eqn:mu_i:def}
\mu_i(t)=\left\{\begin{array}{ll}\lambda_{n,\lfloor (i-1)/2\rfloor+1}(t), & i\in\mathbb{O}, \\
\lambda_{\lfloor (i-1)/2\rfloor+1, n}(t), & i\in\mathbb{E}, \end{array} \right.
\end{equation}
and
\begin{equation} \label{eqn:G_i:def}
G_i=\left\{\begin{array}{ll}G'_{n,\lfloor (i-1)/2\rfloor+1}, & i\in\mathbb{O}, \\
G'_{\lfloor (i-1)/2\rfloor+1, n}, & i\in\mathbb{E}. \end{array} \right.
\end{equation}
The Poisson counters $\{N_i(t)\}_{t\in\mathbb{R}_{\geq 0}}$, $1\leq i\leq m$, determine the rates of jump between the states of the Markov chain in~(\ref{eqn:uvr:sys}). Now, noting that this Markov chain models the sampling instances $\{T_i^\ell\}_{i=0}^\infty$, $1\leq \ell\leq L$, using the jumps that occur in its state $x(t)$, we can control the average sampling frequencies of the sensors through the rates~$\mu_i(t)$, $1\leq i\leq m$. Similar to~\cite{4738725}, we assume that we can control the rates as
\begin{equation} \label{eqn:markovcontrol}
\mu_i(t)=\mu_{i,0}+\sum_{j=1}^m \alpha_{ij} u_j(t),
\end{equation}
and thereby control the average sampling frequencies\footnote{Notice that the number of control inputs in~\eqref{eqn:markovcontrol} does not have to be the same as the number of Poisson counters for the proofs in Subsection~\ref{sec:OptimalSchedulingPolicy} to hold. However, we decided to follow the convention of~\cite{4738725} because we use these results in Sections~\ref{sec:NetworkedEstimationPerformance} and~\ref{sec:NetworkedControlPerformance} to optimally schedule sensors in a networked system in which we can control all the rates $(\mu_i(t))_{i=1}^N$ directly. }.  In~\eqref{eqn:markovcontrol}, $\alpha_{ij}\in\mathbb{R}$, $1\leq i,j\leq N$, are constant parameters that determine the sensitivity of Poisson counters' jump rates with respect to control inputs $u_j$ for $1\leq j\leq N$. Control signals $u_j(t)$, $1\leq j\leq m$, are chosen in order to minimize the cost function
\begin{equation} \label{eqn:uvr:2}
J=\lim_{T\rightarrow\infty} \mathbb{E}\bigg\{\frac{1}{T}\int_{0}^T \sum_{\ell=1}^L \xi_\ell \,e_n^\top x(t)\;\mathrm{d}N_{2\ell}(t)+u(t)^\top u(t) \;\mathrm{d}t \bigg\},
\end{equation}
where $\xi_\ell\in\mathbb{R}_{\geq 0}$, $1\leq \ell\leq L$, are design parameters. Note that the cost function~(\ref{eqn:uvr:2}) consists of two types of terms: $\frac{1}{T}\int_{0}^T e_n^\top x(t)\mathrm{d}N_{2\ell}(t)$ denotes the average frequency of the jumps from $I$ to $S_\ell$ in the Markov chain (i.e., the average sampling frequency of sensor~$\ell$); and $\frac{1}{T}\int_0^T u(t)^\top u(t)\mathrm{d}t$ penalizes the control effort in regulating this frequency. If the latter term is removed, the problem would become ill-posed as the optimal rates $\mu_i(t)$ then is zero and $\mathbb{E}\{\mathrm{d}N_i(t)\}=0$. Consequently, the average sampling frequencies of the sensors vanish.

Considering the identity $\mathbb{E}\{\mathrm{d}N_{2\ell}(t)\}=(\mu_{2\ell,0}+\sum_{j=1}^m \alpha_{2\ell,j} u_j(t))\mathrm{d}t$, we can rewrite the cost function in~(\ref{eqn:uvr:2}) as
\begin{equation}\label{eqn:uvr:3}
\begin{split}
J=\lim_{T\rightarrow\infty} \mathbb{E}\bigg\{\frac{1}{T}\int_{0}^T c^\top x(t) + u(t)^\top S x(t)+u(t)^\top u(t)\;\mathrm{d}t\bigg\},
\end{split}
\end{equation}
where $c=e_n\sum_{\ell=1}^L \xi_\ell\mu_{2\ell,0}$ and $S\in\mathbb{R}^{m\times n}$ is a matrix whose entries are defined as $s_{ji}=\sum_{\ell=1}^L \xi_\ell \alpha_{2\ell,j}$ if $i=n$ and $s_{ji}=0$ otherwise.  In the rest of this paper, we use the notation $S_i$, $1\leq i\leq m$, to denote $i$-th row of matrix $S$. In the next subsection, we find a policy that minimizes~(\ref{eqn:uvr:3}) with respect to the rate control law~(\ref{eqn:markovcontrol}) and subject to the Markov chain dynamics~(\ref{eqn:uvr:sys}). Doing so, we develop an optimal scheduling policy which fairly allocates the network resources (i.e., the sampling instances) between the devices in a sensor network.

\subsection{Optimal Sensor Scheduling} \label{sec:OptimalSchedulingPolicy}
We start by minimizing the finite horizon version of the cost function in~(\ref{eqn:uvr:2}). The proof of the following theorem is a slight generalization of Brockett's result in~\cite{4738725} but follows the same line of reasoning\footnote{The statement makes use of the concept of infinitesimal generators. See~\cite[pp.\,124]{Oksendal2003stochastic} for definition and discussion.}.

\begin{theorem} \label{tho:1} Consider a continuous-time Markov chain evolving on $\mathcal{X}=\{e_1,\dots,e_n\}\subset \mathbb{R}^n$, generated by~(\ref{eqn:uvr:sys}). Let us define matrices $A=\sum_{i=1}^m \mu_{i,0} G_i$ and $B_i=\sum_{j=1}^m \alpha_{ji}G_j$, where for all $1\leq i,j\leq m$, $G_i$ and $\alpha_{ij}$ are introduced in~(\ref{eqn:G_i:def}) and~(\ref{eqn:markovcontrol}), respectively. Assume that, for given $T\in\mathbb{R}_{>0}$ and $c:[0,T]\rightarrow \mathbb{R}^n$, the differential equation
\begin{equation} \label{eqn:diff:eqn:k(t)}
\dot{k}(t)=-c(t)-A^\top k(t)+\frac{1}{4}\sum_{i=1}^m (S_i^\top +B_i^\top k(t))^{.2}; \; k(T)=k_f,
\end{equation}
has a solution on $[0,T]$ such that, for each $(t,x)\in[0,T]\times \mathcal{X}$, the operator
$A-\sum_{i=1}^m \frac{1}{2} (k(t)^\top B_i+S_i) x B_i$ is an infinitesimal generator. Then, the control law
\begin{equation} \label{eqn:tho:1:optimalcontrol}
u_i(t,x)=-\frac{1}{2}\left(k(t)^\top B_i+S_i\right)x(t),\hspace{.1in}1\leq i\leq m,
\end{equation}
minimizes
\begin{equation*} \label{eqn:modcost}
\begin{split}
J=\mathbb{E}&\left\{\frac{1}{T}\int_{0}^T c(t)^\top x(t)+ u(t)^\top Sx(t)+u(t)^\top u(t) \mathrm{d}t +\frac{1}{T}k_f^\top x(T) \right\}.
\end{split}
\end{equation*}
Furthermore, $J=\frac{1}{T}k(0)^\top \mathbb{E}\left\{x(0) \right\}$.
\end{theorem}

\begin{myproof} We follow a similar reasoning as in~\cite{4738725} to calculate the optimal Poisson rates. By adding and subtracting the term $\left.k(t)^\top \mathbb{E}\left\{x(t) \right\}\right|_{0}^T$ from the right hand-side of the scaled cost function $TJ-k_f^\top \mathbb{E}\left\{x(T) \right\}$, we get
\begin{equation} \label{eqn:1}
\begin{split}
TJ-k_f^\top \mathbb{E}\left\{x(T) \right\}&=\mathbb{E}\left\{\int_{0}^Tc(t)^\top x(t)+ u(t)^\top Sx(t)+u(t)^\top u(t) \mathrm{d}t \right\}\\&=-\left.k(t)^\top \mathbb{E}\left\{x(t) \right\}\right|_{0}^T+\left.k(t)^\top \mathbb{E}\left\{x(t) \right\}\right|_{0}^T\\ &\hspace{.93in}+\mathbb{E}\left\{\int_{0}^Tc(t)^\top x(t)+ u(t)^\top Sx(t)+u(t)^\top u(t) \mathrm{d}t \right\}.
\end{split}
\end{equation}
Using the identity
$\left.k(t)^\top \mathbb{E}\left\{x(t) \right\}\right|_{0}^T=\mathbb{E}\left\{\int_{0}^T\mathrm{d}\langle k(t),x(t)\rangle \right\}$ inside~(\ref{eqn:1}), we get
\begin{equation} \label{eqn:2}
\begin{split}
TJ-k_f^\top \mathbb{E}\left\{x(T) \right\}=&-\left.k(t)^\top \mathbb{E}\left\{x(t) \right\}\right|_{0}^T+\mathbb{E}\left\{\int_{0}^T\mathrm{d}\langle k(t),x(t)\rangle \right\}\\&\hspace{.4in}+\mathbb{E}\left\{\int_{0}^Tc(t)^\top x(t)+ u(t)^\top Sx(t)+u(t)^\top u(t) \mathrm{d}t \right\}.
\end{split}
\end{equation}
Using It\^{o}'s Lemma~\cite[p.\,49]{Oksendal2003stochastic}, we know that
$$
\mathrm{d}\langle k(t),x(t)\rangle=\langle \dot{k}(t),x(t)\rangle \mathrm{d}t+\sum_{i=1}^m \langle k(t),G_i x(t) \; \mathrm{d} N_i(t)\rangle,
$$
which transforms~(\ref{eqn:2}) into
\begin{equation*}
\begin{split}
TJ-k_f^\top \mathbb{E}\left\{x(T) \right\}=-\left.k(t)^\top \mathbb{E}\left\{x(t) \right\}\right|_{0}^T&+\mathbb{E}\left\{\int_{0}^T\langle \dot{k}(t),x(t)\rangle \mathrm{d}t+\sum_{i=1}^m \langle k(t),G_i x(t) \;\mathrm{d} N_i(t)\rangle \right\}\\&+\mathbb{E}\left\{\int_{0}^Tc(t)^\top x(t)+ u(t)^\top Sx(t)+u(t)^\top u(t) \mathrm{d}t \right\}.
\end{split}
\end{equation*}
Taking expectation over $x(t)$ and the Poisson processes $\{N_i(t)\}_{t\in\mathbb{R}_{\geq 0}}$, $1\leq i\leq m$, we get
\begin{equation} \label{eqn:3}
\begin{split}
TJ-k_f^\top \mathbb{E}\left\{x(T) \right\}=&-\left.k(t)^\top \mathbb{E}\left\{x(t) \right\}\right|_{0}^T + \int_{0}^T \langle \dot{k}(t)+c(t)+A^\top k(t),p(t)\rangle \mathrm{d}t \\& \hspace{.4in}+ \mathbb{E}\left\{ \int_0^T \hspace{-.03in} u(t)^\top u(t) \hspace{-.03in}+ \hspace{-.03in}\sum_{i=1}^m u_i(t) (S_ix(t) \hspace{-.03in}+\hspace{-.03in}\langle k(t),B_i x(t) \rangle) \mathrm{d}t\right\}\hspace{-.03in},
\end{split}
\end{equation}
where, for $1\leq i\leq m$, $S_i$ is $i$-th row of matrix $S$ and $p(t)=\mathbb{E}\{x(t)\}$. We can rewrite~(\ref{eqn:3}) as
\begin{equation} \label{eqn:cost:positive:term}
\begin{split}
TJ-k_f^\top \mathbb{E}\{x(T) \}=& \int_{0}^T \langle \dot{k}(t)+c(t)+A^\top k(t)-\frac{1}{4}\sum_{i=1}^m (S_i^\top +B_i^\top k(t))^{.2},p(t)\rangle \mathrm{d}t \\& -\left.k(t)^\top \mathbb{E}\left\{x(t) \right\}\right|_{0}^T + \mathbb{E}\left\{\int_0^T \sum_{i=1}^m \left\|u_i(t)+\frac{1}{2}(k(t)^\top B_i+S_i)x(t) \right\|^2 \mathrm{d}t \right\},
\end{split}
\end{equation}
using completion of squares. As there exists a well-defined solution to the differential equation~(\ref{eqn:diff:eqn:k(t)}), the first integral in~(\ref{eqn:cost:positive:term}) vanishes. Hence, the optimal control law is given by~(\ref{eqn:tho:1:optimalcontrol}) since this control law minimizes the last term of~(\ref{eqn:cost:positive:term}). Consequently, equation~(\ref{eqn:cost:positive:term}) gives
$$
TJ=k_f^\top \mathbb{E}\left\{x(T) \right\}-\left.k(t)^\top \mathbb{E}\left\{x(t) \right\}\right|_{0}^T=k(0)^\top \mathbb{E}\left\{x(0) \right\}.
$$
This completes the proof.
\end{myproof}

Notice that for some parameter settings of the cost function, the operator
$A-\sum_{i=1}^m \frac{1}{2} (k(t)^\top B_i+S_i) x B_i$ may not be an infinitesimal generator. A future avenue of research could be to characterize these cases and to present conditions for avoiding them. 

Based on Theorem~\ref{tho:1}, we are able to solve the following infinite-horizon version of the optimal scheduling policy. In the infinite-horizon case, we need to assume that the parameters of the Markov chain and the cost function are time invariant.

\begin{corollary} \label{cor:3} Consider a continuous-time Markov chain evolving on $\mathcal{X}=\{e_1,\dots,e_n\}\subset \mathbb{R}^n$, generated by~(\ref{eqn:uvr:sys}). Let us define matrices $A=\sum_{i=1}^m \mu_{i,0} G_i$ and $B_i=\sum_{j=1}^m \alpha_{ji}G_j$, where for all $1\leq i,j\leq m$, $G_i$ and $\alpha_{ij}$ are introduced in~(\ref{eqn:G_i:def}) and~(\ref{eqn:markovcontrol}), respectively. Assume that, for a given $c\in\mathbb{R}^n$, the nonlinear equation
\begin{equation} \label{eqn:nonlinearequationcor:3}
\matrix{cc}{A^\top & -\mathrm{\textbf{1}} \\ \textbf{1}^\top & 0} \matrix{c}{k_0 \\ \varrho }-
\frac{1}{4}\matrix{c}{\sum_{i=1}^m (S_i^\top +B_i^\top k_0)^{.2} \\ 0 }=\matrix{c}{-c\\ 0},
\end{equation}
has a solution $(k_0,\varrho)\in\mathbb{R}^n\times \mathbb{R}$ such that, for all $x\in \mathcal{X}$, the operator $A-\sum_{i=1}^m \frac{1}{2}(k_0^\top B_i+S_i) x B_i$ is an infinitesimal generator. Then, the control law
\begin{equation} \label{eqn:infop}
u_i(t,x)=-\frac{1}{2}(k_0^\top B_i+S_i)x(t), \hspace{.1in} 1\leq i\leq m,
\end{equation}
minimizes
\begin{equation*}
J=\lim_{T\rightarrow\infty}\mathbb{E}\left\{\frac{1}{T}\int_{0}^T c^\top x(t)+ u(t)^\top Sx(t)+u(t)^\top u(t) \mathrm{d}t \right\}.
\end{equation*}
Furthermore, we have $J=\varrho$.
\end{corollary}

\begin{myproof} Since $x(t)\in \mathcal{X}$ is bounded (because $\|x(t)\|_2\equiv1$ for $t\in\mathbb{R}_{\geq 0}$), we get the identity
\begin{equation} \label{eqn:long_equation}
\begin{split}
\lim_{T\rightarrow\infty}\mathbb{E}&\left\{\frac{1}{T}\int_{0}^T c^\top x(t)+ u(t)^\top Sx(t)+u(t)^\top u(t) \mathrm{d}t \right\}\\ &= \lim_{T\rightarrow\infty}\mathbb{E}\left\{\frac{1}{T}\int_{0}^T c^\top x(t)+ u(t)^\top Sx(t)+u(t)^\top u(t) \mathrm{d}t+\frac{1}{T} k_0^\top x(T)\right\}.
\end{split}
\end{equation}
According to Theorem~\ref{tho:1}, in order to minimize~(\ref{eqn:long_equation}) for any fixed $T\in\mathbb{R}_{>0}$, we have
\begin{equation} \label{eqn:cor2:1}
\dot{k}(t)=-c(t)-A^\top k(t)+\frac{1}{4}\sum_{i=1}^m (S_i^\top +B_i^\top k(t))^{.2},
\end{equation}
with the final condition~$k(T)=k_0$. Defining $q(t)=k(T-t)-k_0-\varrho \textbf{1}t$, we get
\begin{equation*} 
\begin{split}
\dot{q}(t)=& \;-\dot{k}(T-t)-\varrho \textbf{1}\\[-.5em]=&\;A^\top k(T-t)+c-\frac{1}{4}\sum_{i=1}^m (S_i^\top+ B_i^\top k(T-t))^{.2}-\varrho \textbf{1}\\[-.5em]=&\;A^\top (q(t)+k_0+\varrho\textbf{1}t)+c-\varrho \textbf{1} -\frac{1}{4}\sum_{i=1}^m (S_i^\top+B_i^\top (q(t)+k_0+\varrho\textbf{1}t))^{.2}.
\end{split}
\end{equation*}
Note that $A^\top \textbf{1}=0$ and $B_i^\top \textbf{1}=0$, $1\leq i\leq m$, as $A$ and $B_i$ are infinitesimal generators. Hence, 
\begin{equation} \label{eqn:diff:eqn:q(t)}
\begin{split}
\dot{q}(t)=&\;A^\top (q(t)+k_0)+c-\varrho \textbf{1}-\frac{1}{4}\sum_{i=1}^m (S_i^\top+B_i^\top (q(t)+k_0))^{.2}\\[-.5em]=& \;A^\top q(t)-\frac{1}{4}\sum_{i=1}^m (S_i^\top+B_i^\top (q(t)+k_0))^{.2}+\frac{1}{4}\sum_{i=1}^m (S_i^\top+B_i^\top k_0)^{.2}.
\end{split}
\end{equation}
Notice that $q^*=0$ is an equilibrium of~(\ref{eqn:diff:eqn:q(t)}), so $q(t)=0$ for all $t\in[0,T]$ since $q(0)=k(T)-k_0=0$. Therefore, we get $k(t)=k_0 + \varrho \textbf{1}(T-t)$, which results in
$\frac{1}{2}(k(t)^\top B_i+S_i)=\frac{1}{2}(k_0^\top B_i+S_i),$
since $\textbf{1}^\top B_i =0$, $1\leq i\leq m$. As a result, when $T$ goes to infinity, the controller which minimizes~(\ref{eqn:long_equation}) is given by~(\ref{eqn:infop}).
Furthermore, we have
\begin{equation*}
J=\lim_{T\rightarrow\infty} \frac{1}{T} k(0)^\top \mathbb{E}\{x(0)\}=
\lim_{T\rightarrow\infty} \frac{1}{T} \big(k_0 + \varrho \textbf{1}(T-0)\big)^\top \mathbb{E}\{x(0)\} = \varrho\textbf{1}^\top\mathbb{E}\{x(0)\}=\varrho.
\end{equation*}
Finally, notice that the condition $\mathbf{1}^T k_0 = 0$ in the second row of~\eqref{eqn:nonlinearequationcor:3} reduces the number of solutions $k_0$ that satisfy the nonlinear equation in the first row of~\eqref{eqn:nonlinearequationcor:3}. Removing this condition, $k_0+\vartheta\mathbf{1}$ for any $\vartheta\in\mathbb{R}$ is a solution. Notice that all these parallel solutions result in the same control law because $((k_0+\vartheta\mathbf{1})^\top B_i+S_i)=(k_0^\top B_i+S_i)$ following the fact that $\mathbf{1}^\top B_i=0$ for all $1\leq i\leq N$. 
\end{myproof}

Corollary~\ref{cor:3} introduces an optimal scheduling policy to fairly allocate measurement transmissions among sensors according to the cost function in~(\ref{eqn:uvr:2}). By changing the design parameters~$\xi_\ell$, $1\leq \ell\leq L$, we can tune the average sampling frequencies of the subsystems according to their performance requirements. In addition, by adding an extra term to the cost function whenever a new subsystem in introduced or by removing a term whenever a subsystem is detached, we can easily accommodate for dynamic changes in an ad-hoc network. In the remainder of this section, we analyze the asymptotic properties of the optimal scheduling policy in Corollary~\ref{cor:3}.

\subsection{Average Sampling Frequencies}
In this subsection, we study the relationship between the Markov chain parameters and the effective sampling frequencies of the subsystems. Recalling from the problem formulation, $\{T_i^\ell\}_{i=0}^\infty$ denotes the sequence of time instances that the state of the Markov chain in~(\ref{eqn:uvr:1}) jumps from the idle node $I$ to $S_\ell$ and hence, subsystem~$\ell$ is sampled. Mathematically, we define these time instances as
$$
T_0^\ell=\inf\{t\geq 0\;|\; \exists \;\epsilon>0: x(t-\epsilon)=e_n \wedge x(t)=e_\ell\}, 
$$
and
$$
T_{i+1}^\ell=\inf\{t\geq T_i^\ell \;|\; \exists \;\epsilon>0: x(t-\epsilon)=e_n \wedge x(t)=e_\ell\}, \hspace{.1in} i\in\mathbb{Z}_{\geq 0}.
$$
Furthermore, we define the sequence of random variables $\{\Delta_i^\ell\}_{i=0}^\infty$ such that $\Delta_i^\ell=T_{i+1}^\ell-T_i^\ell$ for all $i\in\mathbb{Z}_{\geq 0}$. These random variables denote the time interval between any two successive sampling instances of sensor~$\ell$. We make the assumption that the first and second samples happen within finite time almost surely:

\begin{assumption} \label{asm:1} $\mathbb{P}\{T_0^\ell<\infty\}=1$ and $\mathbb{P}\{T_1^\ell<\infty\}=1$. \end{assumption}

This assumption is not restrictive. Note that it is trivially satisfied if the number of subsystems is finite, the Markov chain is irreducible, and the rates of Poisson processes are finite and uniformly bounded away from zero.

\begin{lemma} \label{lem:1} $\{\Delta_i^\ell\}_{i=0}^\infty$ are identically and independently distributed random variables. \end{lemma}

\begin{myproof} According to the Markov property~\cite[p.\,117]{Oksendal2003stochastic}, we know that, for a given $x(T_i^\ell)$, the trajectory $\{x(t)\,|\,t\geq T_i^\ell\}$ is independent of the history $\{x(t)\,|\,t<T_i^\ell\}$. Noting that $x(T_i^\ell)=e_\ell$ for all $i\geq 1$, gives that $\{\Delta_i^\ell\}_{i=0}^\infty$ are independent random variables. In addition, the Markov chain in~(\ref{eqn:uvr:sys}) and the control law in~(\ref{eqn:infop}) are time invariant. Therefore, the closed-loop Markov chain is also time invariant, and as a result, $\{\Delta_i^\ell\}_{i=0}^\infty$ have equal probability distributions.
\end{myproof}

Now, we are ready to prove that the average sampling frequency of subsystems~$\ell$ is actually equal to $\lim_{T\rightarrow \infty}\frac{1}{T}\int_{0}^T e_n^\top x(t)\mathrm{d}N_{2\ell}(t)$. However, first, we prove the following useful lemma.

\begin{lemma} \label{tho:3} Let the sequence of sampling instances $\{T_i^\ell\}_{i=0}^\infty$ satisfy Assumption~\ref{asm:1}. Then, 
\begin{equation} \label{eqn:tho:1}
\lim_{t\rightarrow \infty} \frac{M_t^\ell}{t}\stackrel{as}{=}\frac{1}{\mathbb{E}\{\Delta_i^\ell\}},
\end{equation}
where $M_t^\ell=\max\left\{i\geq 1\;|\; T_i^\ell\leq t \right\}$ counts the number of jumps prior to any given time $t\in\mathbb{R}_{\geq 0}$ and $x\stackrel{as}{=}y$ means that $\mathbb{P}\{x=y\}=1$.
\end{lemma}

\begin{myproof} This proof follows a similar reasoning as in the proof of Theorem~14 in~\cite{LectureB}. For any given $\vartheta\in\mathbb{Z}_{> 0}$, we have
$$
T_{\vartheta}^\ell=T_0^\ell+\sum_{i=1}^\vartheta T_{i}^\ell-T_{i-1}^\ell=T_0^\ell+\sum_{i=1}^\vartheta \Delta_{i-1}^\ell.
$$
Note that since $\mathbb{P}\{T_0^\ell<\infty\}=1$ according to Assumption~\ref{asm:1}, we get
$\lim_{\vartheta\rightarrow \infty} T_0^\ell/\vartheta\stackrel{as}{=}0$. Therefore, we have
\begin{equation*}
\lim_{\vartheta\rightarrow\infty} \frac{T_{\vartheta}^\ell}{\vartheta}\stackrel{as}{=}
\lim_{\vartheta\rightarrow\infty} \frac{1}{\vartheta} \sum_{i=1}^\vartheta \Delta_{i-1}^\ell .
\end{equation*}
Notice that $\mathbb{P}\{\Delta_0^\ell=T_1^\ell-T_0^\ell<\infty\}=1$ according to Assumption~\ref{asm:1}. Therefore, using Lemma~\ref{lem:1}, we get $\mathbb{P}\{\Delta_i^\ell<\infty\}=1$ for all $i\in\mathbb{Z}_{\geq 0}$. Consequently, $\mathbb{E}\{|\Delta_i^\ell|\}= \mathbb{E}\{\Delta_i^\ell\}<\infty$. Now, using the strong law of large numbers~\cite{laha1979probability}, we get
\begin{equation} \label{eqn:proof}
\lim_{\vartheta\rightarrow\infty} \frac{T_{\vartheta}^\ell}{\vartheta}\stackrel{as}{=}
\mathbb{E}\{\Delta_{i}^\ell\}.
\end{equation}
For any $t\in\mathbb{R}_{\geq 0}$, we have $T_{M_t^\ell}^\ell\leq t < T_{M_t^\ell+1}^\ell.$
Therefore, we get $T_{M_t^\ell}^\ell/M_t^\ell\leq t/M_t^\ell<T_{M_t^\ell+1}^\ell/M_t^\ell.$
Notice that $\lim_{t\rightarrow \infty} M_t^\ell\stackrel{as}{=}\infty$ since, as proved earlier, $\mathbb{P}\{\Delta_i^\ell<\infty\}=1$ for all $i\in\mathbb{Z}_{\geq 0}$. Using~(\ref{eqn:proof}), we get
$$
\lim_{t\rightarrow \infty}\frac{T_{M_t^\ell}^\ell}{M_t^\ell}\stackrel{as}{=}\mathbb{E}\{\Delta_i^\ell\},
\hspace{.3in} \lim_{t\rightarrow \infty}\frac{T_{M_t^\ell+1}^\ell}{M_t^\ell}\stackrel{as}{=}\mathbb{E}\{\Delta_i^\ell\},
$$
which results in~(\ref{eqn:tho:1}).
\end{myproof}

We now state our main result concerning the average sampling frequency of the sensors denoted by
$$
f_\ell=\lim_{T\rightarrow \infty}\mathbb{E}\left\{\frac{1}{T}\int_{0}^T e_n^\top x(t)\;\mathrm{d}N_{2\ell}(t)\right\}, \hspace{.1in } 1\leq\ell\leq L.
$$

\begin{theorem} \label{cor:2} Let the sequence of sampling instances $\{T_i^\ell\}_{i=0}^\infty$ satisfy Assumption~\ref{asm:1}. If $\lim_{t\rightarrow \infty} p(t)$ exists, the average sampling frequency of sensor~$\ell$ is equal to
\begin{equation*}
\begin{split}
f_\ell&=\frac{1}{\mathbb{E}\{\Delta_i^\ell\}}=\left(\mu_{2\ell,0} -\frac{1}{2}\sum_{j=1}^m \alpha_{2\ell, j}(k_0^\top B_j+S_j)e_n\right) e_n^\top \lim_{t\rightarrow \infty} p(t),
\end{split}
\end{equation*}
where $p(t)=\mathbb{E}\{x(t)\}$ can be computed by
\begin{equation} \label{eqn:tho:ptdynamics}
\dot{p}(t)=\left(A-\frac{1}{2}\sum_{i=1}^m B_i \Lambda(k_0^\top B_i+S_i) \right) p(t), \; p(0)=\mathbb{E}\left\{x(0)\right\},
\end{equation}
with notation $\Lambda(k_0^\top B_i+S_i)=\diag((k_0^\top B_i+S_i)e_1,\dots,(k_0^\top B_i+S_i)e_n)$.
\end{theorem}

\begin{myproof} The proof of equality $f_\ell=1/\mathbb{E}\{\Delta_i^\ell\}$ directly follows from applying Lemma~\ref{tho:3} in conjunction with that $M_T^\ell=\int_{0}^T e_n^\top x(t) \mathrm{d}N_{2\ell}(t)$. Now, we can compute $p(t)$ using
\begin{equation} \label{eqn:ptevolution}
\dot{p}(t)=Ap(t)+\mathbb{E}\left\{\sum_{i=1}^m u_i(t,x(t)) B_i x(t)\right\}, \; p(0)=\mathbb{E}\left\{x(0)\right\}
\end{equation}
Substituting~(\ref{eqn:infop}) inside~(\ref{eqn:ptevolution}), we get
\begin{equation*}
\begin{split}
\dot{p}(t)&=Ap(t) - \frac{1}{2}\mathbb{E}
\left\{ \sum_{i=1}^m (k_0^\top B_i+S_i) x(t) B_i x(t) \right\}
\\& = Ap(t) - \frac{1}{2}\mathbb{E}
\left\{ \sum_{i=1}^m (k_0^\top B_i+S_i) \matrix{c}{ x_1(t) \\  \vdots  \\ x_n(t) }   B_i  \matrix{c}{ x_1(t) \\  \vdots  \\  x_n(t) }\right\} \\& = Ap(t) - \frac{1}{2}\mathbb{E}
\left\{ \sum_{i=1}^m B_i  \matrix{c}{x_1(t)\sum_{j=1}^n (k_0^\top B_i+ S_i)e_j x_j(t) \\ \vdots \\ x_n(t)\sum_{j=1}^n (k_0^\top B_i+S_i)e_j x_j(t)} \right\} .
\end{split}
\end{equation*}
Note that $x_{\zeta}(t) \sum_{j=1}^n (k_0^\top B_i+S_i)e_j x_j(t)=(k_0^\top B_i+S_i)e_{\zeta} x_{\zeta}(t)$ for $1\leq \zeta\leq n$, since $x(t)\in\mathcal{X}$ is a unit vector in $\mathbb{R}^n$. Therefore, we get~(\ref{eqn:tho:ptdynamics}). Now, noticing that $p(t)$ converges exponentially to a nonzero steady-state value as time goes to infinity (because otherwise $\lim_{t\rightarrow \infty} p(t)$ does not exist), we can expand the expression for the average sampling frequencies of the sensors as
\begin{equation} \label{eqn:simplefl}
\begin{split}
f_\ell&=\lim_{T\rightarrow \infty}
\mathbb{E}\left\{\frac{1}{T}\int_{0}^T e_n^\top x(t)\left(\mu_{2\ell,0}+\sum_{j=1}^m \alpha_{2\ell, j}u_{j}\right) \mathrm{d}t\right\}\\&=\lim_{T\rightarrow \infty}\mathbb{E}\left\{
\frac{1}{T}\int_{0}^T e_n^\top x(t)\left(\mu_{2\ell,0}-\frac{1}{2}\sum_{j=1}^m \alpha_{2\ell, j}(k_0^\top B_j+S_j)x(t)\right) \mathrm{d}t\right\}
\\&=\lim_{T\rightarrow \infty}
\frac{1}{T}\int_{0}^T e_n^\top p(t)\left(\mu_{2\ell,0}-\frac{1}{2}\sum_{j=1}^m \alpha_{2\ell, j}(k_0^\top B_j+S_j)e_n\right) \mathrm{d}t
\\&=\left(\mu_{2\ell,0}-\frac{1}{2}\sum_{j=1}^m \alpha_{2\ell, j}(k_0^\top B_j+S_j)e_n\right) e_n^\top \lim_{t\rightarrow \infty} p(t),
\end{split}
\end{equation}
where the third equality follows again from the fact that $x(t)\in\mathcal{X}$ is a unit vector.
\end{myproof}

Theorem~\ref{cor:2} allows us to calculate the average sampling frequencies of the subsystems. We use these average sampling frequencies to bound the closed-loop performance of the networked system when the proposed optimal scheduling policy is implemented.

\section{Applications to Networked Estimation} \label{sec:NetworkedEstimationPerformance}
In this section, we study networked estimation based on the proposed stochastic scheduling policy. Let us start by presenting the system model and the estimator. As a starting point, we introduce a networked system that is composed of scalar decoupled subsystems. In Subsections~\ref{eqn:arbitrarydim} and~\ref{subsec:Kalmanfiltering}, we generalize some of the results to decoupled higher-order subsystems.

\subsection{System Model and Estimator}
Consider the networked system illustrated in Figure~\ref{figureNCS_new}, where subsystem~$\ell$, $1\leq \ell\leq L$, is a scalar stochastic system described by
\begin{equation} \label{eqn:sys}
\mathrm{d}z_\ell(t)=-\gamma_\ell z_\ell(t) \;\mathrm{d}t+\sigma_\ell \;\mathrm{d}w_\ell(t); \; z_\ell(0)=0,
\end{equation}
with given model parameters $\gamma_\ell,\sigma_\ell\in\mathbb{R}_{\geq 0}$. Note that all subsystems are stable. The stochastic processes $\{w_\ell(t)\}_{t\in\mathbb{R}_{\geq 0}}$, $1\leq \ell \leq L$, are statistically independent Wiener processes with zero mean. Estimator~$\ell$ receives state measurements $\{y_i^\ell\}_{i=0}^\infty$ at time instances~$\{T_i^\ell\}_{i=0}^\infty$, such that
\begin{equation} \label{eqn:output}
y_i^\ell=z_\ell(T_i^\ell)+n_i^\ell; \hspace{.2in}  \forall i\in\mathbb{Z}_{\geq 0},
\end{equation}
where $\{n_i^\ell\}_{i=0}^\infty$ denotes measurement noise sequence, which is composed of independently and identically distributed Gaussian random variables with zero mean and specified standard deviation~$\eta_\ell$. Let each subsystem adopt a simple estimator of the form
\begin{equation} \label{eqn:estimation}
\frac{\mathrm{d}}{\mathrm{d}t}\hat{z}_\ell(t)=-\gamma_\ell \hat{z}_\ell(t);\hspace{.2in} \hat{z}_\ell(T_i^\ell)=y_i^\ell,
\end{equation}
for $t\in[T_i^\ell,T_{i+1}^\ell)$. We define the estimation error $e_\ell(t)=z_\ell(t)-\hat{z}_\ell(t)$. Estimator~$\ell$ only has access to the state measurements of subsystem~$\ell$ at specific time instances~$\{T_i^\ell\}_{i=0}^\infty$ but is supposed to reconstruct the signal at any time $t\in\mathbb{R}_{\geq 0}$.  Notice that this estimator is not optimal. In Subsection~\ref{subsec:Kalmanfiltering}, we will consider estimators based on Kalman filtering instead.

\subsection{Performance Analysis: Scalar Subsystems}
In this subsection, we present an upper bound for the performance of the introduced networked estimator. First, we prove the following simple lemma.

\begin{lemma} \label{lem:2} Let the function $g:\mathbb{R}_{\geq 0}\rightarrow \mathbb{R}$ be defined as $g(t)=c_1e^{-2\gamma t}+\frac{c_2}{2\gamma}(1-e^{-2\gamma t})$ with given scalars $c_1,c_2\in\mathbb{R}$ and $\gamma\in\mathbb{R}_{> 0}$ such that $2\gamma c_1 \leq c_2$. Then,
\\ (a) $g$ is a non-decreasing function on its domain;
\\ (b) $g$ is a concave function on its domain.
\end{lemma}

\begin{myproof}
For part~(a), note that if $2\gamma c_1 \leq c_2$, the function $g(t)$ is continuously differentiable and $\textrm{d}g(t)/\textrm{d}t =-(2\gamma c_1-c_2) e^{-2\gamma t} \geq 0$ for all $t\in\mathbb{R}_{\geq 0}$. Hence, $g(t)$ is a non-decreasing function on its domain (since it is continuous). On the other hand, for part~(b), note that if $2\gamma c_1 \leq c_2$, the function $g(t)$ is double continuously differentiable and $\textrm{d}^2g(t)/\textrm{d}t^2 =2\gamma(2\gamma c_1-c_2) e^{-2\gamma t} \leq 0 $ for all $t\in\mathbb{R}_{\geq 0}$. Therefore, $g(t)$ is a concave function on its domain.
\end{myproof}

The following theorem presents upper bounds for the estimation error variance for the cases where the measurement noise is small or large, respectively.

\begin{theorem} \label{tho:estimation:1} Assume that subsystem~$\ell$, $1\leq \ell \leq L$, is described by~(\ref{eqn:sys}) and let the sequence of sampling instances $\{T_i^\ell\}_{i=0}^\infty$ satisfy Assumption~\ref{asm:1}. Then, if $\eta_\ell \leq \sqrt{1/(2\gamma_\ell)}\sigma_\ell$, the estimation error variance is bounded by
\begin{equation} \label{eqn:upperbound:refined:scalar_error}
\begin{split}
\mathbb{E}\{e_\ell^2(t)\} \leq \eta_\ell^2 e^{-2\gamma_\ell/f_\ell}+\frac{\sigma_\ell^2}{2\gamma_\ell}\left(1 - e^{-2\gamma_\ell/f_\ell}\right),
\end{split}
\end{equation}
otherwise, if $\eta_\ell > \sqrt{1/(2\gamma_\ell)}\sigma_\ell$,
\begin{equation} \label{eqn:scalar_error}
\begin{split}
\mathbb{E}\{e_\ell^2(t)\} \leq \eta_\ell^2 +\frac{\sigma_\ell^2}{2\gamma_\ell}\left(1 - e^{-2\gamma_\ell/f_\ell}\right).
\end{split}
\end{equation}
\end{theorem}

\begin{myproof} Using It\^{o}'s Lemma~\cite[p.\,49]{Oksendal2003stochastic}, for all $t\in[T_i^\ell,T_{i+1}^\ell)$, we get
\begin{equation*}
\begin{split}
\mathrm{d}e_\ell(t)&=\left(-\frac{\mathrm{d}}{\mathrm{d}t}\hat{z}_\ell(t)-\gamma_\ell e_\ell(t) -\gamma_\ell \hat{z}_\ell(t) \right)\mathrm{d}t+\sigma_\ell \mathrm{d}w_\ell(t)
=-\gamma_\ell e_\ell(t)\mathrm{d}t+\sigma_\ell \mathrm{d}w_\ell(t),
\end{split}
\end{equation*}
with the initial condition $e_\ell(T_i^\ell)=-n_i^\ell$. First, let us consider the case where $\eta_\ell \leq \sqrt{1/(2\gamma_\ell)}\sigma_\ell$. Again, using It\^{o}'s Lemma, we get
\begin{equation*}
\begin{split}
\mathrm{d}(e_\ell^2(t))&=(-2\gamma_\ell \,e_\ell^2(t)+\sigma_\ell^2)\mathrm{d}t+2e_\ell(t)\sigma_\ell \mathrm{d}w_\ell(t),
\end{split}
\end{equation*}
and as a result
$$
\frac{\mathrm{d}}{\mathrm{d}t}\mathbb{E}\{e_\ell^2(t)\,|\,\Delta_i^\ell\}=-2\gamma_\ell \,\mathbb{E}\{e_\ell^2(t)\,|\,\Delta_i^\ell\}+\sigma_\ell^2,
$$
where $\mathbb{E}\{e_\ell^2(T_i^\ell)\,|\,\Delta_i^\ell\}=\eta_\ell^2$.
Hence, for all $t\in[T_i^\ell,T_{i+1}^\ell)$, we have
\begin{equation*}
\begin{split}
\mathbb{E}\{e_\ell^2(t)\,|\,\Delta_i^\ell \}= \eta_\ell^2 e^{-2\gamma_\ell (t-T_i^\ell)}+\frac{\sigma_\ell^2}{2\gamma_\ell}\left(1 - e^{-2\gamma_\ell (t-T_i^\ell)}\right).
\end{split}
\end{equation*}
Now, using Lemma~\ref{lem:2}~(a), it is easy to see that
\begin{equation*}
\begin{split}
\mathbb{E}\{e_\ell^2(t)\,|\,\Delta_i^\ell \}&\leq \eta_\ell^2 e^{-2\gamma_\ell \Delta_i^\ell}+\frac{\sigma_\ell^2}{2\gamma_\ell}\left(1 - e^{-2\gamma_\ell \Delta_i^\ell}\right).
\end{split}
\end{equation*}
Note that
\begin{equation} \label{eqn:long_equation_1}
\begin{split}
\mathbb{E}\{e_\ell^2(t)\} & = \mathbb{E}\{\mathbb{E}\{e_\ell^2(t)\,|\,\Delta_i^\ell\}\}
 \leq \mathbb{E}\bigg\{\eta_\ell^2 e^{-2\gamma_\ell \Delta_i^\ell}+\frac{\sigma_\ell^2}{2\gamma_\ell}\left(1 - e^{-2\gamma_\ell \Delta_i^\ell}\right)\bigg\}.
\end{split}
\end{equation}
By using Lemma~\ref{lem:2}~(b) along with Jensen's Inequality~\cite[p.\,320]{Oksendal2003stochastic}, we can transform~(\ref{eqn:long_equation_1}) into~(\ref{eqn:upperbound:refined:scalar_error}). For the case where $\eta_\ell > \sqrt{1/(2\gamma_\ell)}\sigma_\ell$, we can similarly derive the upper bound
\begin{equation*}
\begin{split}
\mathbb{E}\{e_\ell^2(t)\,|\,\Delta_i^\ell \}\leq \eta_\ell^2 +\frac{\sigma_\ell^2}{2\gamma_\ell}\left(1 - e^{-2\gamma_\ell \Delta_i^\ell}\right),
\end{split}
\end{equation*}
which results in~(\ref{eqn:scalar_error}), again using Jensen's Inequality.
\end{myproof}

Note that the upper bound~(\ref{eqn:upperbound:refined:scalar_error}) is tighter than~(\ref{eqn:scalar_error}) when the equality~$\eta_\ell =\sqrt{1/(2\gamma_\ell)}\sigma_\ell$ holds. In the next two subsections, we generalize these results to higher-order subsystems. 

\subsection{Performance Analysis: Higher-Order Subsystems with Noisy State Measurement}
\label{eqn:arbitrarydim}
Let us assume that subsystem~$\ell$, $1\leq \ell\leq L$, is described by
\begin{equation} \label{eqn:generalsubsystem}
\mathrm{d}z_\ell(t)=A_\ell z_\ell(t)\mathrm{d}t+ H_\ell \mathrm{d}w_\ell(t); \hspace{.2in} z_\ell(0)=0,
\end{equation}
where $z_\ell(t)\in\mathbb{R}^{d_\ell}$ is its state with $d_\ell\in\mathbb{Z}_{\geq 1}$ and $A_\ell$ is its model matrix satisfying $\overline{\lambda}(A_\ell+A_\ell^\top)<0$ where $\overline{\lambda}(\cdot)$ denotes the largest eigenvalue of a matrix. In addition, $\{w_\ell(t)\}_{t\in\mathbb{R}_{\geq 0}}$, $1\leq \ell \leq L$, is a tuple of statistically independent Wiener processes with zero mean. Estimator~$\ell$ receives noisy state-measurements~$\{y_i^\ell\}_{i=0}^\infty$ at time instances~$\{T_i^\ell\}_{i=0}^\infty$, such that
\begin{equation} \label{eqn:output:generalsubssytem}
y_i^\ell=z_\ell(T_i^\ell)+n_i^\ell; \hspace{.1in}  \forall i\in\mathbb{Z}_{\geq 0},
\end{equation}
where $\{n_i^\ell\}_{i=0}^\infty$ denotes the measurement noise and is composed of independently and identically distributed Gaussian random variables with $\mathbb{E}\{n_i^\ell\}=0$ and $\mathbb{E}\{n_i^\ell n_i^{\ell\top}\}=R_\ell$. We define the estimation error as $e_\ell(t)=z_\ell(t)-\hat{z}_\ell(t)$, where for all $t\in[T_i^\ell,T_{i+1}^\ell)$, the state estimate $\hat{z}_\ell(t)$ is derived by
$$
\frac{\mathrm{d}}{\mathrm{d}t}\hat{z}_\ell(t)=A_\ell\hat{z}_\ell(t);\; \hat{z}_\ell(T_i^\ell)=y_i^\ell,
$$
The next theorem presents an upper bound for the variance of this estimation error. For scalar subsystems, the introduced upper bound in~(\ref{eqn:upperbound:nonscalar_error}) is equivalent to the upper bound in~(\ref{eqn:scalar_error}).

\begin{theorem} \label{tho:estimation:2} Assume that subsystem~$\ell$, $1\leq \ell \leq L$, is described by~(\ref{eqn:generalsubsystem}) and let the sequence of sampling instances $\{T_i^\ell\}_{i=0}^\infty$ satisfy Assumption~\ref{asm:1}. Then, the estimation error variance is bounded by
\begin{equation} \label{eqn:upperbound:nonscalar_error}
\begin{split}
\mathbb{E}\{\|e_\ell(t)\|^2\}\leq \trace(R_\ell)+\frac{\trace(H^\top H)}{|\overline{\lambda}(A_\ell+A_\ell^\top)|}
\left(1-e^{\overline{\lambda}(A_\ell+A_\ell^\top)/f_\ell}\right).
\end{split}
\end{equation}
\end{theorem}

\begin{myproof}
Using It\^{o}'s Lemma, for all $t\in[T_i^\ell,T_{i+1}^\ell)$, we get
\begin{equation*}
\begin{split}
\mathrm{d}\|e_\ell(t)\|^2=&\;e_\ell(t)^\top (A_\ell+A_\ell^\top)e_\ell(t)\mathrm{d}t+\trace(H^\top H)\mathrm{d}t+e_\ell(t)^\top H \mathrm{d}w_\ell(t)+\mathrm{d}w_\ell(t)^\top H^\top e_\ell(t),
\end{split}
\end{equation*}
and as a result
\begin{equation*}
\begin{split}
\frac{\mathrm{d}}{\mathrm{d}t}\mathbb{E}\{\|e_\ell(t)\|^2\,|\,\Delta_i^\ell\}
&=\trace(H^\top H)+\mathbb{E}\{e_\ell(t)^\top (A_\ell+A_\ell^\top)e_\ell(t)\,|\,\Delta_i^\ell\}
\\&\leq \trace(H^\top H) +\overline{\lambda}(A_\ell+A_\ell^\top) \mathbb{E}\{\|e_\ell(t)\|^2\,|\,\Delta_i^\ell\},
\end{split}
\end{equation*}
with the initial condition $\mathbb{E}\{\|e_\ell(T_i^\ell)\|^2\}=\trace(R_\ell)$. Now, using the Comparison Lemma~\cite[p.102]{khalil2002nonlinear}, we get
\begin{equation*}
\begin{split}
\mathbb{E}\{\|e_\ell(t)\|^2\,|\,\Delta_i^\ell\}\leq \trace(R_\ell) e^{\overline{\lambda}(A_\ell+A_\ell^\top)t}+
\frac{\trace(H^\top H)}{|\overline{\lambda}(A_\ell+A_\ell^\top)|} \left(1-e^{\overline{\lambda}(A_\ell+A_\ell^\top)t}\right),
\end{split}
\end{equation*}
for $t\in[T_i^\ell,T_{i+1}^\ell)$. Using Lemma~\ref{lem:2} and Jensen's Inequality, we get~(\ref{eqn:upperbound:nonscalar_error}).
\end{myproof}

It is possible to refine the upper bound~(\ref{eqn:upperbound:nonscalar_error}) for the case where $\trace(R_\ell)  \leq 1/(|\overline{\lambda}(A_\ell+A_\ell^\top)|) \linebreak[4] \trace(H^\top H)$, following a similar argument as in the proof of Theorem~\ref{tho:estimation:1}.

\subsection{Performance Analysis: Higher-Order Subsystems with Noisy Output Measurement}
\label{subsec:Kalmanfiltering}
In this subsection, we assume that estimator~$\ell$, $1\leq \ell \leq L$, receives noisy output measurements~$\{y_i^\ell\}_{i=0}^\infty$ at time instances~$\{T_i^\ell\}_{i=0}^\infty$, such that
\begin{equation} \label{eqn:output:generalsubsystem:output}
y_i^\ell=C_\ell z_\ell(T_i^\ell)+n_i^\ell; \hspace{.1in}  \forall i\in\mathbb{Z}_{\geq 0},
\end{equation}
where $C_\ell\in\mathbb{R}^{p_\ell\times d_\ell}$ (for a given output vector dimension $p_\ell\in\mathbb{Z}_{\geq 1}$ such that $p_\ell\leq d_\ell$) and the measurement noise $\{n_i^\ell\}_{i=0}^\infty$ is a sequence of independently and identically distributed Gaussian random variables with $\mathbb{E}\{n_i^\ell\}=0$ and $\mathbb{E}\{n_i^\ell n_i^{\ell\,\top}\}=R_\ell$. For any sequence of sampling instances~$\{T_i^\ell\}_{i=0}^\infty$, we can discretize the stochastic continuous-time system in~(\ref{eqn:generalsubsystem}) as
$$
z_\ell[i+1]=F_\ell[i]z_\ell[i]+G_\ell[i]w_\ell[i],
$$
where $z_\ell[i]=z(T_i^\ell)$, $F_\ell[i]=e^{A(T_{i+1}^\ell-T_i^\ell)}$, and the sequence $\{G_\ell[i]\}_{i=0}^\infty$ is chosen such that
$$
G_\ell[i]G_\ell[i]^\top=\int_0^{T_{i+1}^\ell-T_i^\ell} e^{A\tau}HH^\top e^{A^\top \tau}\mathrm{d}\tau,
\hspace{.2in} \forall i\in\mathbb{Z}_{\geq 0}.
$$
In addition, $\{w_\ell[i]\}_{i=0}^\infty$ is a sequence of independently and identically distributed Gaussian random variables with zero mean and unity variance. It is evident that $y_\ell[i]=C_\ell z_\ell[i]+n_i^\ell$. We run a discrete-time Kalman filter over these output measurements to calculate the state estimates~$\{\hat{z}_\ell[i]\}_{i=0}^\infty$ with error covariance matrix $P_\ell[i]=\mathbb{E}\{(z_\ell[i]-\hat{z}_\ell[i])(z_\ell[i]-\hat{z}_\ell[i])^\top\}$. For inter-sample times $t\in[T_i^\ell,T_{i+1}^\ell)$, we use a simple prediction filter
\begin{equation} \label{eqn:esimator:Kalmanbased}
\frac{\mathrm{d}}{\mathrm{d}t}\hat{z}_\ell(t)=A_\ell\hat{z}_\ell(t);
\;\hat{z}_\ell(T_i^\ell)=\hat{z}_\ell[i].
\end{equation}
Let us define the estimation error as $e_\ell(t)=z_\ell(t)-\hat{z}_\ell(t)$. The next theorem present an upper bound for the estimation error variance.

\begin{theorem} \label{tho:estimation:3} Assume that subsystem~$\ell$, $1\leq \ell \leq L$, is described by~(\ref{eqn:generalsubsystem}). Then, the estimator given by~(\ref{eqn:esimator:Kalmanbased}) is an optimal mean square error estimator and for any fixed sequence of sampling instances $\{T_i^\ell\}_{i=0}^\infty$, the estimation error is upper-bounded by 
\begin{equation} \label{eqn:error:Kalman:Fixedtimes}
\begin{split}
\mathbb{E}\{&\|e_\ell(t)\|^2\,|\,\Delta_i^\ell\} \leq \mathrm{trace}(P_\ell[i]) +
\frac{\trace(H^\top H)}{|\overline{\lambda}(A_\ell+A_\ell^\top)|} \left(1-e^{\overline{\lambda}(A_\ell+A_\ell^\top)(T_{i+1}^\ell-T_i^\ell)}\right).
\end{split}
\end{equation}
\end{theorem}

\begin{myproof} First, note that for $t\in[T_i^\ell,T_{i+1}^\ell)$, the estimator
$$
\frac{\mathrm{d}}{\mathrm{d}t}\hat{z}_\ell(t)=A_\ell\hat{z}_\ell(t);
\;\hat{z}_\ell(T_i^\ell)=\mathbb\{z_\ell(T_i^\ell)\,|\,y_1^\ell,\dots,y_i^\ell\},
$$
is an optimal mean square error estimator. This is in fact true since the estimator~$\ell$ has not received any new information over $[T_i^\ell,t]$ and it should simply predict the state using the best available estimation $\mathbb\{z_\ell(T_i^\ell)\,|\,y_1^\ell,\dots,y_i^\ell\}$. Now, recalling from~\cite{kailath2000linear}, we know that $\mathbb\{z_\ell(T_i^\ell)\,|\,y_1^\ell,\dots,y_i^\ell\}= \mathbb\{z_\ell[i]\,|\,y_1^\ell,\dots,y_i^\ell\}= \hat{z}_\ell[i]$. This completes the first part of the proof. For the rest, note that following a similar reasoning as in the proof of Theorem~\ref{tho:estimation:2}, for all $t\in[T_i^\ell,T_{i+1}^\ell)$, we get
\begin{equation*}
\begin{split}
\frac{\mathrm{d}}{\mathrm{d}t}\mathbb{E}\{\|e_\ell(t)\|^2\;|\;&\Delta_i^\ell\}\leq \trace(H^\top H) +\overline{\lambda}(A_\ell+A_\ell^\top)
\mathbb{E}\{\|e_\ell(t)\|^2\,|\,\Delta_i^\ell\},
\end{split}
\end{equation*}
with the initial condition
$\mathbb{E}\{\|e_\ell(T_i^\ell)\|^2\}=\mathbb{E}\{(z_\ell[i]-\hat{z}_\ell[i])^\top
(z_\ell[i]-\hat{z}_\ell[i])\}=\mathrm{trace}(P_\ell[i]),$
which results in~(\ref{eqn:error:Kalman:Fixedtimes}) again using the Comparison Lemma.
\end{myproof}

Note that the upper bound~(\ref{eqn:error:Kalman:Fixedtimes}) is conditioned on the sampling intervals. Unfortunately, it is difficult to calculate~$\mathbb{E}\{\mathrm{trace}(P_\ell[i])\}$ as a function of average sampling frequencies, which makes it hard to eliminate the conditional expectation. However, for the case where $p_\ell=n_\ell$, the upper bound~(\ref{eqn:upperbound:nonscalar_error}) would also hold for the estimator in~(\ref{eqn:esimator:Kalmanbased}). This is indeed true because~(\ref{eqn:esimator:Kalmanbased}) is an optimal mean square error estimator.

\section{Applications to Networked Control} \label{sec:NetworkedControlPerformance}
In this section, we study networked control as an application of the proposed stochastic scheduling policy. Let us start by presenting the system model and the control law. We first present the results for impulsive controllers in Subsection~\ref{subsec:impulsive}. However, in Subsections~\ref{subsec:3} and~\ref{subsec:4}, we generalize these results to pulse and exponential controllers.

\subsection{System Model and Controller}
Consider the stochastic control system
\begin{equation} \label{eqn:sys:control}
\mathrm{d}z_\ell(t)=(-\gamma_\ell z_\ell(t)+v_\ell(t)) \;\mathrm{d}t+\sigma_\ell \;\mathrm{d}w_\ell(t); \hspace{.1in} z_\ell(0)=z^0_\ell,
\end{equation}
where $z_\ell(t)\in\mathbb{R}$ and $v_\ell(t)\in\mathbb{R}$, $1\leq \ell\leq L$, are the state and control input of subsystem~$\ell$. We assume that each subsystem is in feedback interconnection with a subcontroller governed by the control law
\begin{equation} \label{eqn:control}
v_\ell(t)=-\sum_{i=0}^\infty y_i^\ell f(t-T^\ell_i),
\end{equation}
where $y_i^\ell=z(T_i^\ell)+n_i$ for all $i\in\mathbb{Z}_{\geq 0}$ and $f:\mathbb{R}\rightarrow \mathbb{R}\cup\{\pm\infty\}$ is chosen appropriately to yield a causal controller (i.e., $f(t)=0$ for all $t<0$). For instance, using $f(\cdot)=\delta(\cdot)$, where $\delta(\cdot)$ is the impulse function (see~\cite[p.\,1]{izrail1964generalized}), results in an impulsive controller, which simply resets the state of its corresponding subsystem to a neighborhood of the origin characterized by the amplitude of the measurement noise whenever a new measurement is received. Without loss of generality, we assume that $z^0_\ell=0$ because the influence of the initial condition is only visible until the first sampling instance $T_0^\ell$, which is guaranteed to happen in a finite time thanks to Assumption~\ref{asm:1}.

\subsection{Performance Analysis: Impulsive Controllers} \label{subsec:impulsive}
In this subsection, we present an upper bound for the closed-loop performance of subsystems described in~(\ref{eqn:sys:control}) and controlled by an impulsive controller. In this case, for all $t\in[T_i^\ell,T_{i+1}^\ell)$, the closed-loop subsystem~$\ell$ is governed by
$$
\mathrm{d}z_\ell(t)=-\gamma_\ell z_\ell(t) \;\mathrm{d}t+\sigma_\ell \;\mathrm{d}w_\ell(t); \; z_\ell(T_i^\ell)=-n_i^\ell.
$$
The next theorem presents an upper bound for the performance of this closed-loop system which corresponds to the estimation error upper bound presented in Theorem~\ref{tho:estimation:1}.

\begin{theorem} \label{tho:control:1} Assume that subsystem~$\ell$, $1\leq \ell \leq L$, is described by~(\ref{eqn:sys:control}) and let the sequence of sampling instances $\{T_i^\ell\}_{i=0}^\infty$ satisfy Assumption~\ref{asm:1}. Then, if $\eta_\ell \leq \sqrt{1/(2\gamma_\ell)}\sigma_\ell$, the closed-loop performance of subsystem~$\ell$ is bounded by
\begin{equation} \label{eqn:result:NS:1}
\begin{split}
\mathbb{E}\left\{z_\ell^2(t) \right\} &\leq  \eta_\ell^2e^{-2\gamma_\ell/f_\ell}+\frac{\sigma_\ell^2}{2\gamma_\ell}
\left(1-e^{-2\gamma_\ell/f_\ell}\right).
\end{split}
\end{equation}
otherwise,
\begin{equation} \label{eqn:highnoise_impulsive}
\begin{split}
\mathbb{E}\left\{z_\ell^2(t) \right\}\leq \eta_\ell^2
+\frac{\sigma_\ell^2}{2\gamma_\ell}\left(1-e^{-2\gamma_\ell/f_\ell}\right).
\end{split}
\end{equation}
\end{theorem}

\begin{myproof} Similar to the proof of Theorem~\ref{tho:estimation:1}. See~\cite{FarokhiCompleteManuscript2012} for details.
\end{myproof}

Note that the closed-loop performance, measured as the variance of the plant state, is upper bounded by the plant and measurement noise variance. In the next two subsections, we generalize this result to pulse and exponential controllers.

\subsection{Performance Analysis: Pulse Controllers}  \label{subsec:3}
In this subsection, we use a narrow pulse function to approximate the behavior of the impulse function. Let us pick a constant $\rho \in\mathbb{R}_{>0}$. For $t\in[T_i^\ell,T_{i+1}^\ell)$, we use the control law
\begin{equation} \label{eqn:almostcontrol}
v_\ell(t)=\left\{\begin{array}{ll} -y_i^\ell\gamma_\ell e^{-\gamma_\ell\rho}/(1-e^{-\gamma_\ell\rho}), & T_i^\ell \leq t\leq T_i^\ell+\rho, \\ 0, & T_i^\ell+\rho < t\leq T_{i+1}^\ell, \end{array}\right.
\end{equation}
whenever $T_i^\ell+\rho \leq T_{i+1}^\ell$, and 
$$
v_\ell(t)=-y_i^\ell\gamma_\ell e^{-\gamma_\ell\rho}/(1-e^{-\gamma_\ell\rho}),\hspace{.2in} T_i^\ell \leq t\leq T_{i+1}^\ell,
$$
otherwise. This controller converges to the impulsive controller as $\rho$ tends to zero.

\begin{theorem} \label{tho:control:2} Assume that subsystem~$\ell$, $1\leq \ell \leq L$, is described by~(\ref{eqn:sys:control}) and let the sequence of sampling instances $\{T_i^\ell\}_{i=0}^\infty$ satisfy Assumption~\ref{asm:1}. Then, the closed-loop performance of subsystem~$\ell$ is bounded by
\begin{equation} \label{eqn:6.1}
\begin{split}
\mathbb{E}\left\{z_\ell^2(t) \right\} &\leq \left[\frac{\sigma_\ell^2}{2\gamma_\ell}\left(1-e^{-2\gamma_\ell/f_\ell}\right)+ \eta_\ell^2 e^{-2\gamma_\ell \rho} \right] \frac{1}{1-\mathbb{P}\{\Delta_i^\ell < \rho\}}.
\end{split}
\end{equation}
\end{theorem}

\begin{myproof} To simplify the calculations, we introduce the change of variable
$z'_\ell(t)=z_\ell(t)+\zeta_\ell(t)$
for all $t\in[T_i^\ell,T_{i+1}^\ell)$, where
\begin{equation*}
\begin{split}
\zeta_\ell(t)=\left\{ \begin{array}{ll} -z_\ell(T_i^\ell)\frac{e^{-\gamma_\ell (t-T_i^\ell)}-e^{-\gamma_\ell \rho}}{1-e^{-\gamma_\ell \rho}}+n_i^\ell \frac{e^{-\gamma_\ell \rho}}{1-e^{-\gamma_\ell \rho}}\left(1-e^{-\gamma_\ell (t-T_i^\ell)}\right), & t\in[T_i^\ell,T_i^\ell+\rho), \\ n_i^\ell e^{-\gamma_\ell (t-T_i^\ell)}, & t\in[T_i^\ell+\rho,T_{i+1}^\ell). \end{array} \right.
\end{split}
\end{equation*}
Now, using It\^{o}'s Lemma~\cite[p.\,49]{Oksendal2003stochastic}, we get 
$$
\mathrm{d}z'_\ell(t)=-\gamma_\ell z'_\ell(t)\mathrm{d}t+\sigma_\ell \mathrm{d}w_\ell(t); \; z'_\ell(T_i^\ell)=0.
$$
Hence, for all $t\in[T_i^\ell,T_{i+1}^\ell)$, we get 
\begin{align*}
\mathbb{E}\left\{z_\ell^2(t) \,|\, \Delta_i^\ell \right\}
& = \mathbb{E}\left\{z_\ell^{\prime 2}(t)+\zeta_\ell^2(t) \,|\, \Delta_i^\ell \right\}
\\&\leq \frac{\sigma_\ell^2}{2\gamma_\ell}\left(1-e^{-2\gamma_\ell \Delta_i^\ell}\right)+\zeta_\ell^2(t),
\end{align*}
where the first equality is due to the fact that $\mathbb{E}\{\zeta_\ell(t)z'_\ell(t)\}=0$ because the random process $\{w_\ell(t)\}_{t\in(T_i^\ell,T_{i+1}^\ell)}$ is independent of $\zeta_\ell(t)$ and $\mathbb{E}\{z'_\ell(t)\}=0$. As a result 
\begin{equation} \label{eqn:5.1}
\begin{split}
\mathbb{E}\left\{z_\ell^2(t) \right\}
\leq\mathbb{E}\left\{\frac{\sigma_\ell^2}{2\gamma_\ell}\left(1-e^{-2\gamma_\ell \Delta_i^\ell}\right) \right\} + \eta_\ell^2 e^{-2\gamma_\ell \rho} +\mathbb{E}\{z_\ell^2(T_i^\ell)\}\mathbb{P}\{\Delta_i^\ell < \rho\}.
\end{split}
\end{equation}
Using Lemma~\ref{lem:2}~(b) and Jensen's Inequality, we can simplify~(\ref{eqn:5.1}) as
\begin{equation} \label{eqn:6.1}
\begin{split}
\mathbb{E}\left\{z_\ell^2(t) \right\}
\leq \frac{\sigma_\ell^2}{2\gamma_\ell}\left(1-e^{-2\gamma_\ell/f_\ell}\right)+ \eta_\ell^2 e^{-2\gamma_\ell \rho}+ \mathbb{E}\{z_\ell^2(T_i^\ell)\}\mathbb{P}\{\Delta_i^\ell < \rho\}.
\end{split}
\end{equation}
Note that by evaluating~\eqref{eqn:6.1} as $t$ goes to $T_{i+1}^\ell$, we can extract a difference equations for the closed-loop performance (i.e., an algebraic equation that relates $\mathbb{E}\{z_\ell^2(T_{i+1}^\ell)\}$ to $\mathbb{E}\{z_\ell^2(T_i^\ell)\}$ for all $i$). By solving this difference equation and substituting the solution into~\eqref{eqn:6.1}, we get  
\begin{equation*}
\begin{split}
\mathbb{E}\left\{z_\ell^2(t) \right\} &\leq \sum_{k=0}^i \left[\frac{\sigma_\ell^2}{2\gamma_\ell}\left(1-e^{-2\gamma_\ell/f_\ell}\right) + \eta_\ell^2 e^{-2\gamma_\ell \rho} \right] \left(\mathbb{P}\{\Delta_i^\ell < \rho\}\right)^k,
\end{split}
\end{equation*}
for all $t\in[T_i^\ell,T_{i+1}^\ell)$, and as a result 
\begin{equation*}
\begin{split}
\mathbb{E}\left\{z_\ell^2(t) \right\} &\leq \sum_{k=0}^\infty \left[\frac{\sigma_\ell^2}{2\gamma_\ell}\left(1-e^{-2\gamma_\ell/f_\ell}\right) + \eta_\ell^2 e^{-2\gamma_\ell \rho} \right] \left(\mathbb{P}\{\Delta_i^\ell < \rho\}\right)^k
\\ &= \left[\frac{\sigma_\ell^2}{2\gamma_\ell}\left(1-e^{-2\gamma_\ell/f_\ell}\right)+ \eta_\ell^2 e^{-2\gamma_\ell \rho} \right] \frac{1}{1-\mathbb{P}\{\Delta_i^\ell < \rho\}}.
\end{split}
\end{equation*}
This concludes the proof.
\end{myproof}

Note that if $\rho$ tends to zero in~(\ref{eqn:6.1}), we would recover the same upper bound as in the case of the impulsive controller~(\ref{eqn:highnoise_impulsive}). This is true since $\lim_{\rho\rightarrow 0}\mathbb{P}\{\Delta_i^\ell < \rho\}=0$ assuming that the probability distribution of hitting-times of the underlying Markov chain is atom-less at the origin, which is a reasonable assumption when the Poisson jump rates are finite. 

\subsection{Performance Analysis: Exponential Controllers}  \label{subsec:4}
In this subsection, we use an exponential function to approximate the impulse function. Let us pick a constant $\theta \in\mathbb{R}_{>0} \setminus \{\gamma_\ell\}$. For all $t\in[T_i^\ell,T_{i+1}^\ell)$, we use the control law 
\begin{equation} \label{eqn:almostcontrol2}
v_\ell(t)=(\gamma_\ell-\theta)y_i^\ell e^{-\theta (t-T_i^\ell)}.
\end{equation}
This controller converges to the impulsive controller as $\theta$ goes to infinity.

\begin{theorem} \label{tho:control:3} Assume that subsystem~$\ell$, $1\leq \ell \leq L$, is described by~(\ref{eqn:sys:control}) and let the sequence of sampling instances $\{T_i^\ell\}_{i=0}^\infty$ satisfy Assumption~\ref{asm:1}. Then, the closed-loop performance of subsystem~$\ell$ is bounded by 
\begin{equation} \label{eqn:result:NS:2}
\begin{split}
\mathbb{E}\left\{z_\ell^2(t) \right\} \leq &\left[\eta_\ell^2+
\frac{\sigma_\ell^2}{2\gamma_\ell}\left(1-e^{-2\gamma_\ell/f_\ell}\right) \right]
\frac{1}{1-\mathbb{E}\{e^{-2\theta\Delta_i^\ell}\}}.
\end{split}
\end{equation}
\end{theorem}

\begin{myproof}
Using the same argument as in the proof of Theorem~\ref{tho:control:2}, we obtain
\begin{equation*}
\begin{split}
\mathbb{E}\left\{z_\ell^2(t) \,|\, \Delta_i^\ell \right\}\leq  \eta_\ell^2&+ \frac{\sigma_\ell^2}{2\gamma_\ell}\left(1-e^{-2\gamma_\ell \Delta_i^\ell}\right)+\mathbb{E}\{z_\ell^2(T_i^\ell)\}e^{-2\theta\Delta_i^\ell},
\end{split}
\end{equation*}
for all $t\in[T_i^\ell,T_{i+1}^\ell)$, and as a result
\begin{equation*}
\begin{split}
\mathbb{E}\left\{z_\ell^2(t) \right\}\leq\eta_\ell^2&+\mathbb{E}\left\{\frac{\sigma_\ell^2}{2\gamma_\ell}\left(1-e^{-2\gamma_\ell \Delta_i^\ell}\right) \right\}+\mathbb{E}\{z_\ell^2(T_i^\ell)\} \mathbb{E}\{e^{-2\theta\Delta_i^\ell}\}.
\end{split}
\end{equation*}
Similar to the proof of Theorem~\ref{tho:control:2}, we can simplify this expression into~(\ref{eqn:result:NS:2}) using Lemma~\ref{lem:2}~(b) and Jensen's Inequality.
\end{myproof}

Note that if $\theta$ goes to infinity, we would recover the same upper bound as in the case of the impulsive controller since $\lim_{\theta\rightarrow +\infty}\mathbb{E}\{e^{-2\theta\Delta_i^\ell}\}=0$ assuming that the probability distribution of hitting-times of the Markov chain is atom-less at the origin. Exponential shape of the control signal is common in biological systems such as in neurological control system~\cite{126254}. 

\section{Numerical Example} \label{sec:NumericalExample}
\begin{figure}[t]
\centering
$$
\begin{array}{cc}
\includegraphics[width=.45\linewidth]{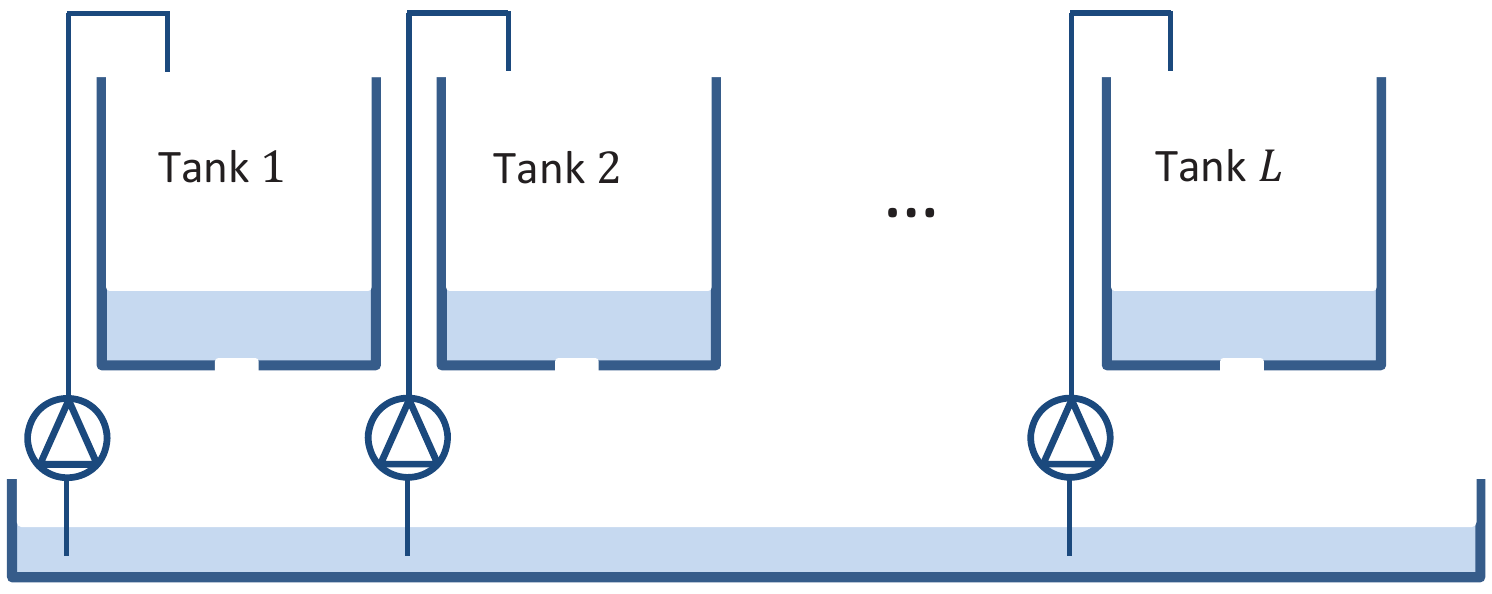} &
\vspace{-.1in}
\includegraphics[width=.32\linewidth]{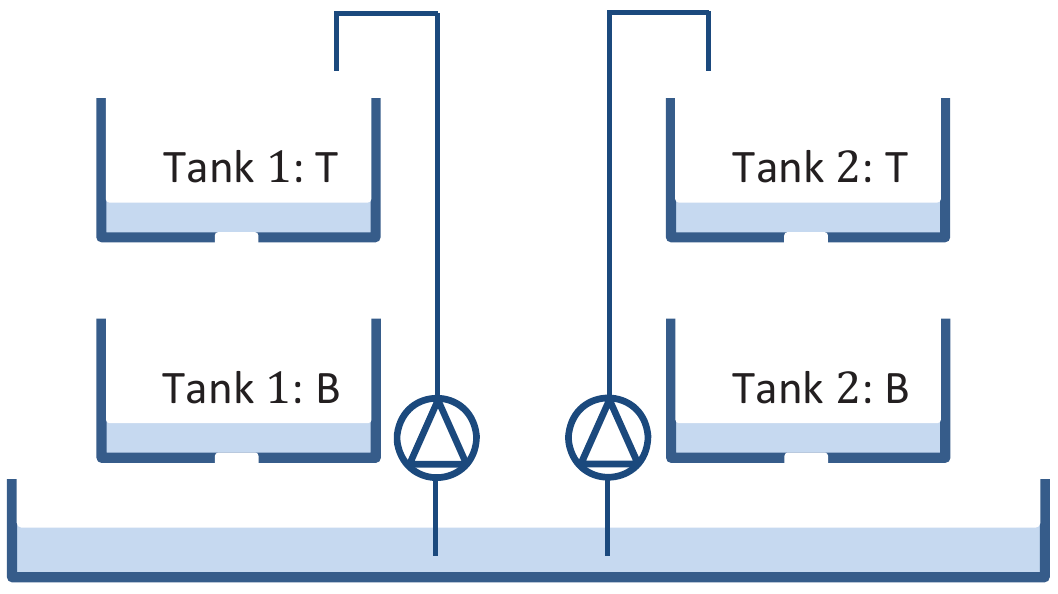}
\vspace{-.1in}
\end{array}
$$
\caption{\label{figurewatertank} An example of a networked system composed of decoupled scalar subsystems (left) and non-scalar subsystems (right). }
\end{figure}

In this section, we demonstrate the developed results on a networked system composed of $L$ decoupled water tanks illustrated in Figure~\ref{figurewatertank}~(left), where each tank is linearized about its stationary water level $h_\ell$ as
$$
\mathrm{d}z_\ell(t)=-\frac{a_\ell}{a'_\ell}\sqrt{\frac{g}{2h_\ell}}z_\ell(t)\mathrm{d}t+\mathrm{d}w_\ell(t); \; z_\ell(0)=z^0_\ell.
$$
In this model, $a'_\ell$ is the cross-section of water tank~$\ell$, $a_\ell$ is the cross-section of its outlet hole, and $g$ is the acceleration of gravity. Furthermore, $z_\ell(t)$ and $v_\ell(t)$ denote the deviation of the tank's water level from its stationary point and the control input, respectively. Let the initial condition $z_\ell^0=0$ as we assume that the tank's water level start at its stationary level. However, due to factors such as input flow fluctuations, the water level drifts away from zero. In the next subsection, we start by numerically demonstrating the estimation results for $L=2$ water tanks.

\subsection{Estimation: Scalar Subsystem} \label{subsec:numerical:scalar:estimation}
Let us fix the parameters $a'_1=a'_2=1.00\,\mathrm{m}^2$, $a_1=0.20\,\mathrm{m}^2$, $a_2=0.10\,\mathrm{m}^2$, $g=9.80\,\mathrm{m}/\mathrm{s}^2$, $h_1=0.40\,\mathrm{m}$, and $h_2=0.54\,\mathrm{m}$. For these physical parameters, the water tanks can be described by~(\ref{eqn:sys}) with $\gamma_1=0.7$, $\gamma_2=0.3$, and $\sigma_1=\sigma_2=1.0$. We sample these subsystems using the Markov chain in~(\ref{eqn:uvr:sys}) with $m=2L=4$. We assume that $\mu_i(t)=\mu_{i,0}+u_i(t)$ for all $1\leq i\leq 4$, where $\mu_{2\ell,0}=1$ and $\mu_{2\ell-1,0}=10$ for $\ell=1,2$. We are interested in finding $u_i(t)$, $1\leq i\leq 4$, in order to minimize the cost function
\begin{equation*} \label{eqn:examplecost}
\begin{split}
J&=\lim_{T\rightarrow \infty}\mathbb{E}\left\{\frac{1}{T}\int_{0}^T 0.5e_3^\top x(t)\mathrm{d}N_2+0.1e_3^\top x(t)\mathrm{d}N_4+u(t)^\top u(t) \mathrm{d}t \right\}.
\end{split}
\end{equation*}
Using Corollary~\ref{cor:3}, we get
\begin{equation*}
\begin{split}
\matrix{c}{u_1(t,x)\\u_2(t,x)\\u_3(t,x)\\u_4(t,x)}=\matrix{ccc}{-0.0228 & 0 &  0 \\0 &  0 & -0.2272
\\0 & -0.0228 &  0\\0 &  0 & -0.0272}x(t).
\end{split}
\end{equation*}

\begin{figure}[t]
\centering
$$
\begin{array}{cc}
\includegraphics[width=.4\linewidth]{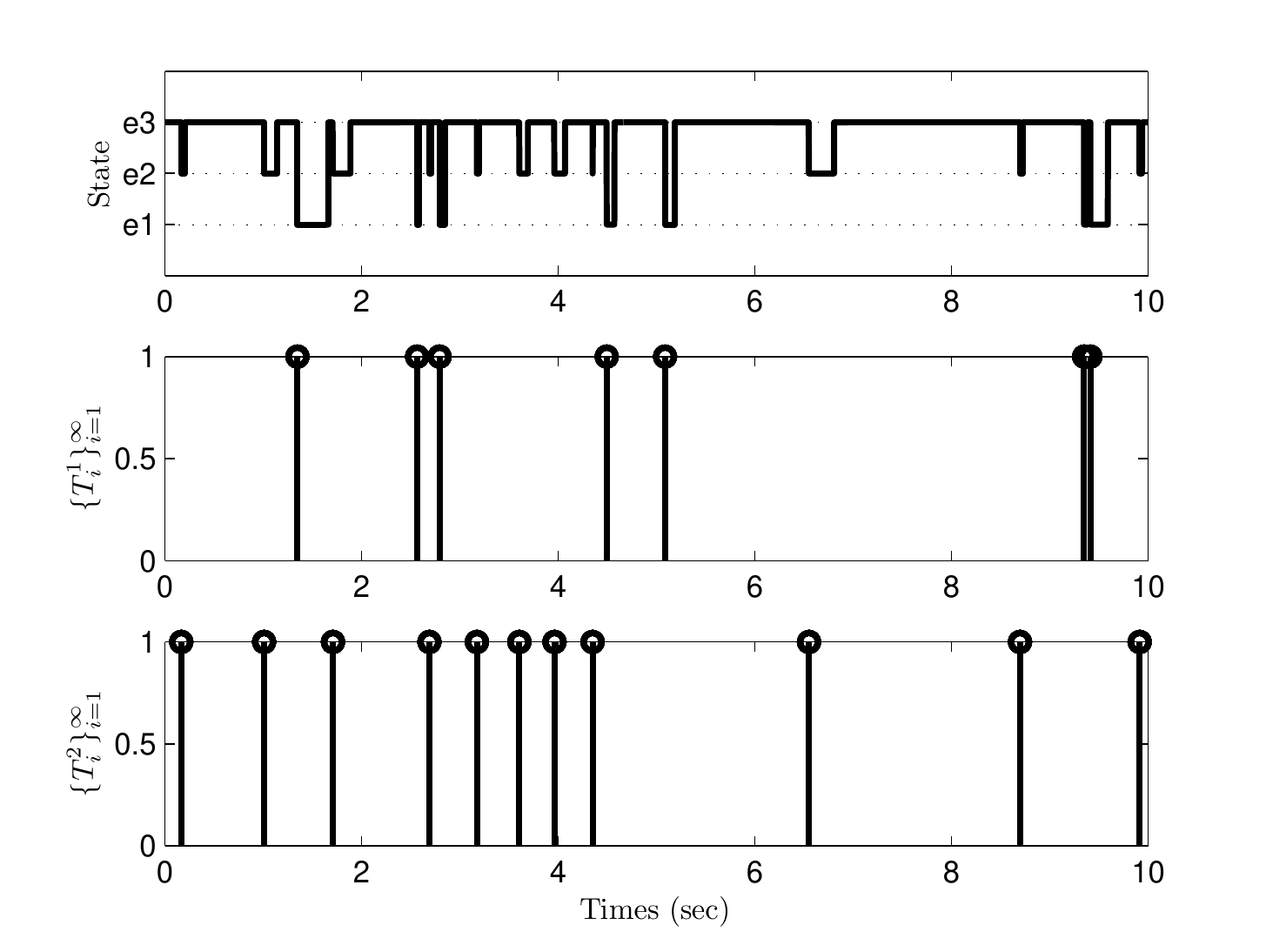} &
\vspace{-.1in}
\includegraphics[width=.4\linewidth]{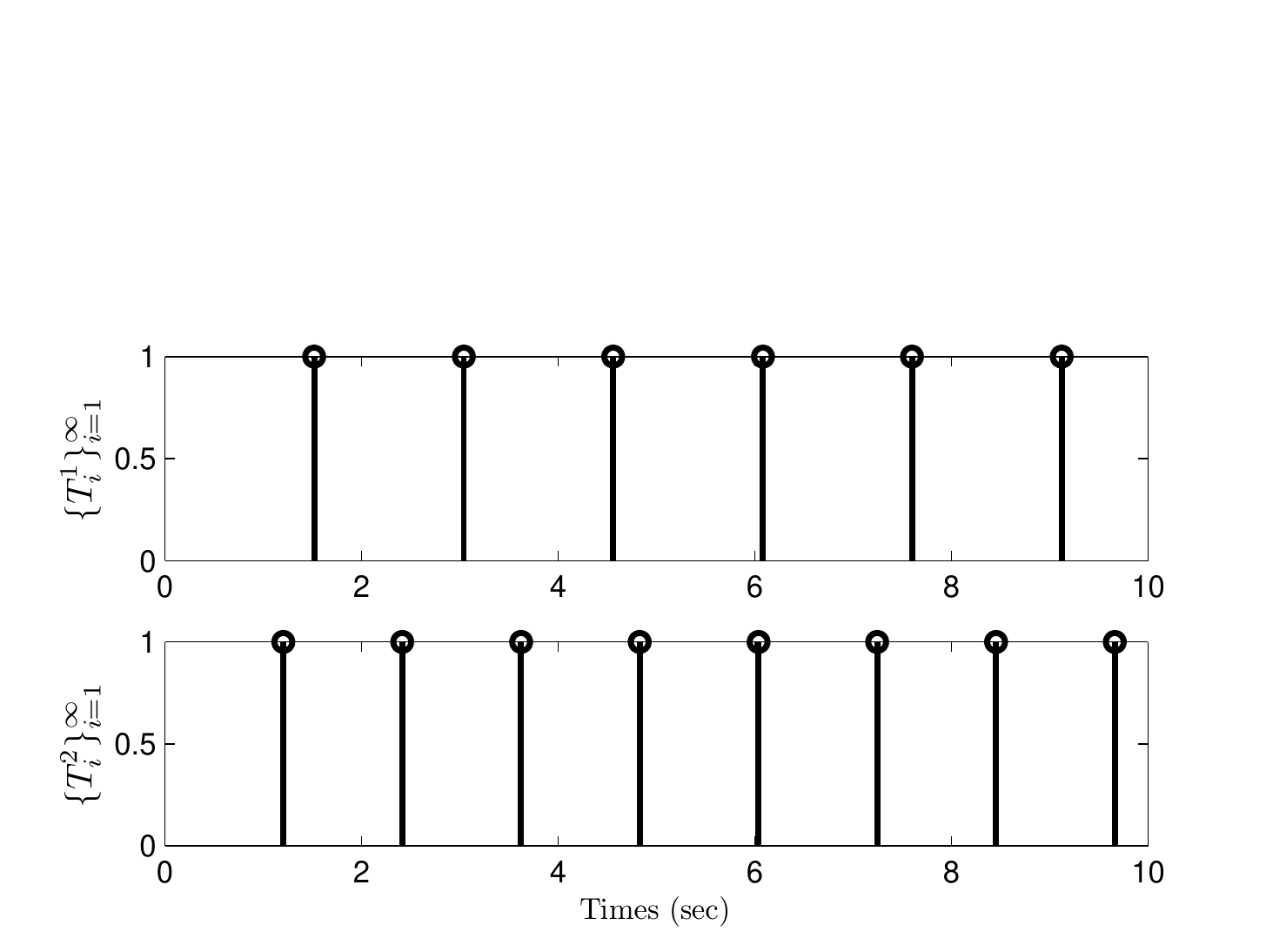} \\
\vspace{-.1in}
\includegraphics[width=.4\linewidth]{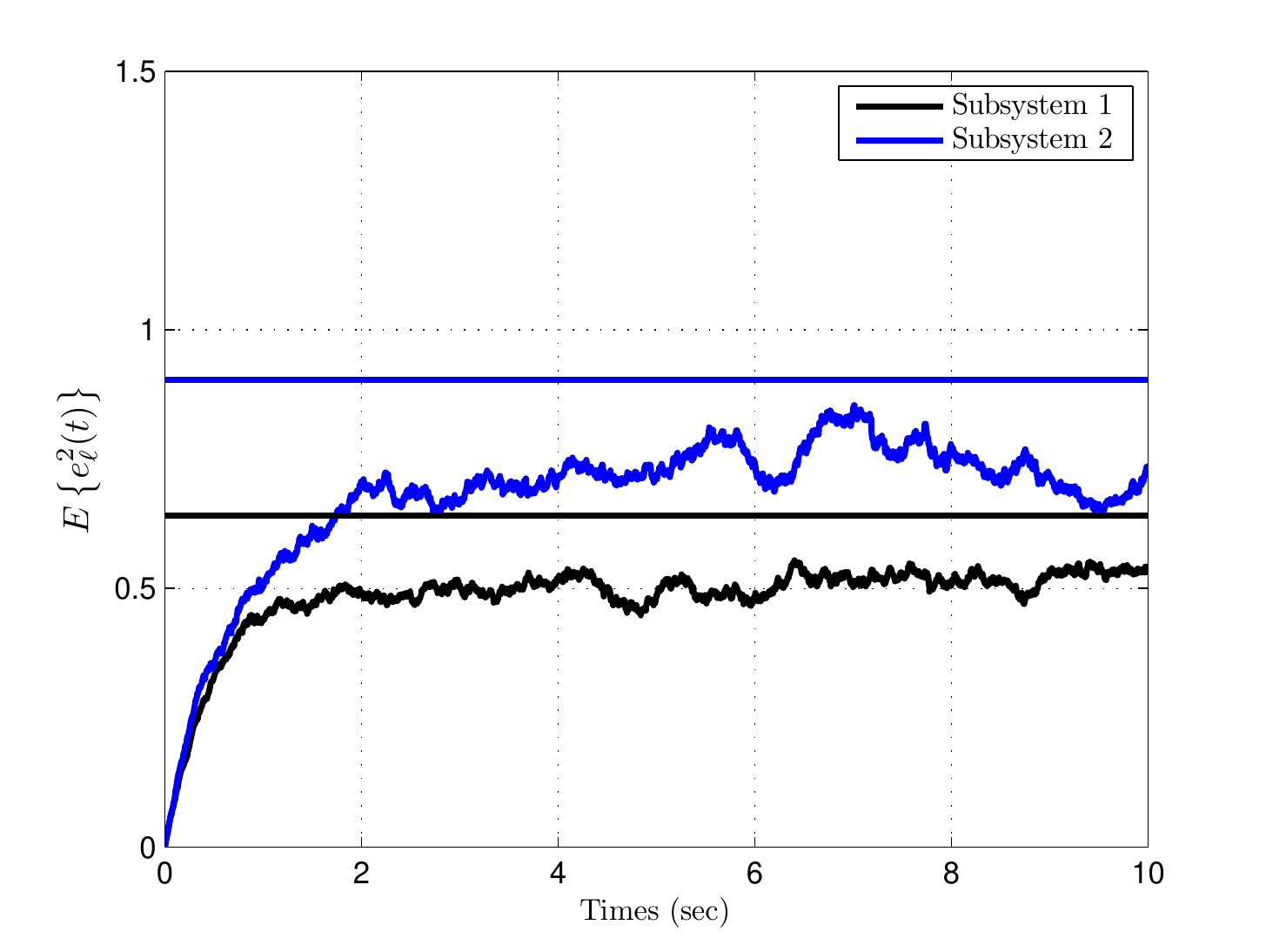} &
\vspace{-.1in}
\includegraphics[width=.4\linewidth]{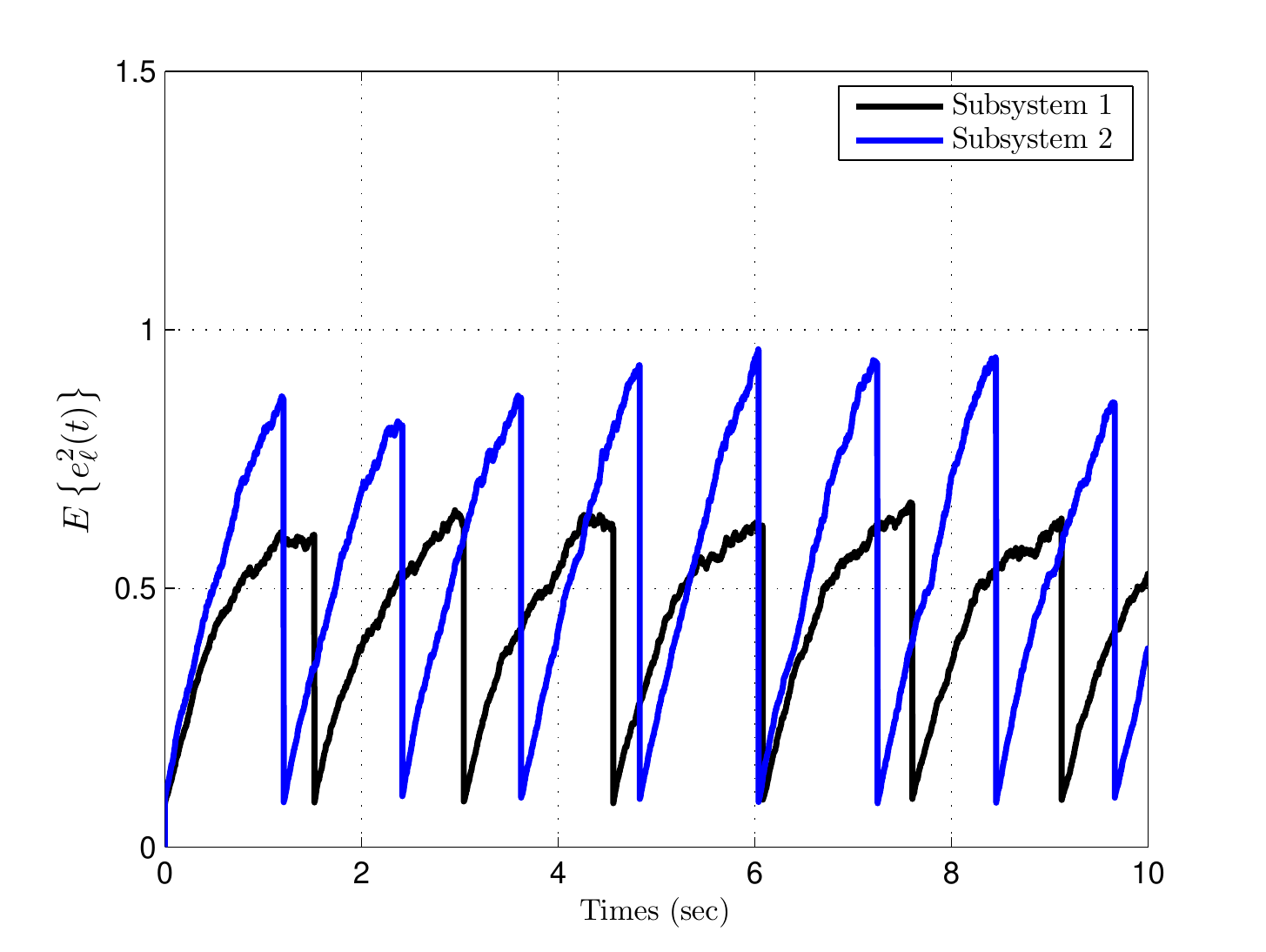}
\vspace{-.1in}
\end{array}
$$
\caption{\label{figureSamplingInstance} An example of the state of the continuous-time Markov chain used in the optimal scheduling policy and its corresponding sampling instances for both subsystems (upper-left). Sampling instances for both subsystems when using a periodic scheduling policy (upper-right). Estimation error $\mathbb{E}\{e_\ell^2(t)\}$ for 1000~Monte Carlo simulations when using the optimal sampling policy (lower-left) and the periodic sampling policy (lower-right). }
\end{figure}

Figure~\ref{figureSamplingInstance}~(upper-left) illustrates an example of the continuous-time Markov chain state $x(t)$ and the sampling instances $\{T_i^\ell\}_{i=0}^\infty$ of subsystems~$\ell=1,2$.  Using~(\ref{eqn:simplefl}), we get the average sampling frequencies $f_1=0.66$ and $f_2=0.83$. Figure~\ref{figureSamplingInstance}~(upper-right) shows the sampling instances when using a periodic scheduling policy with the same sampling frequencies as the average sampling frequencies of the optimal scheduling policy. Note that the optimal scheduling policy allocates the sampling instances according to the jumps between the states of the Markov chain.

We can tune the average sampling frequencies of the subsystems by changing the design parameters $\xi_\ell$, $1\leq \ell \leq L$. Table~\ref{table1} illustrates the average sampling frequencies of the subsystems versus different choices of the design parameters~$\xi_\ell$, $1\leq \ell \leq L$. It is evident that when increasing (decreasing)~$\xi_\ell$ for a given $\ell$, the average sampling frequency of subsystem~$\ell$ decreases (increases).

\begin{table}
\caption{\label{table1} Example of average sampling frequencies. }
\centering
\begin{tabular}{| c | c || c | c |}
  \hline
  $\xi_1$ & $\xi_2$ & $f_1$ & $f_2$ \\
  \hline \hline
  0.1 & 0.1 & 0.8040 & 0.8040 \\
  0.5 & 0.1 & 0.6577 & 0.8279 \\
  1.0 & 0.1 & 0.4656 & 0.8559 \\
  2.0 & 0.1 & 0.0451 & 0.9045 \\
  \hline
\end{tabular}
\end{table}

Let us assume that estimator~$\ell$ has access to state measurements of subsystem~$\ell$ according to~(\ref{eqn:output}) with measurement noise variance $\eta_\ell=0.3$ for $\ell=1,2$. Figure~\ref{figureSamplingInstance}~(lower-left) illustrates the estimation error variance $\mathbb{E}\{e_\ell^2(t)\}$ for 1000~Monte Carlo simulations when using the optimal scheduling policy. The horizontal lines represent the theoretical upper bounds derived in Theorem~\ref{tho:estimation:1}; i.e., $
\mathbb{E}\{e_1^2(t)\} \leq 0.64$ and $\mathbb{E}\{e_2^2(t)\} \leq 0.90$. Note that the approximations of the estimation error variances would eventually converge to the exact expectation value as the number of simulations goes to infinity, and that the theoretical bounds are relatively close. Figure~\ref{figureSamplingInstance}~(lower-right) illustrates the estimation error variance $\mathbb{E}\{e_\ell^2(t)\}$ for 1000~Monte Carlo simulations when using the periodic scheduling policy that is portrayed in Figure~\ref{figureSamplingInstance}~(upper-right). Note that the saw-tooth behavior is due to the fact that the sampling instances are fixed in advance and they are identical for each Monte Carlo simulation.

\begin{figure}[t]
\centering
\includegraphics[width=.4\linewidth]{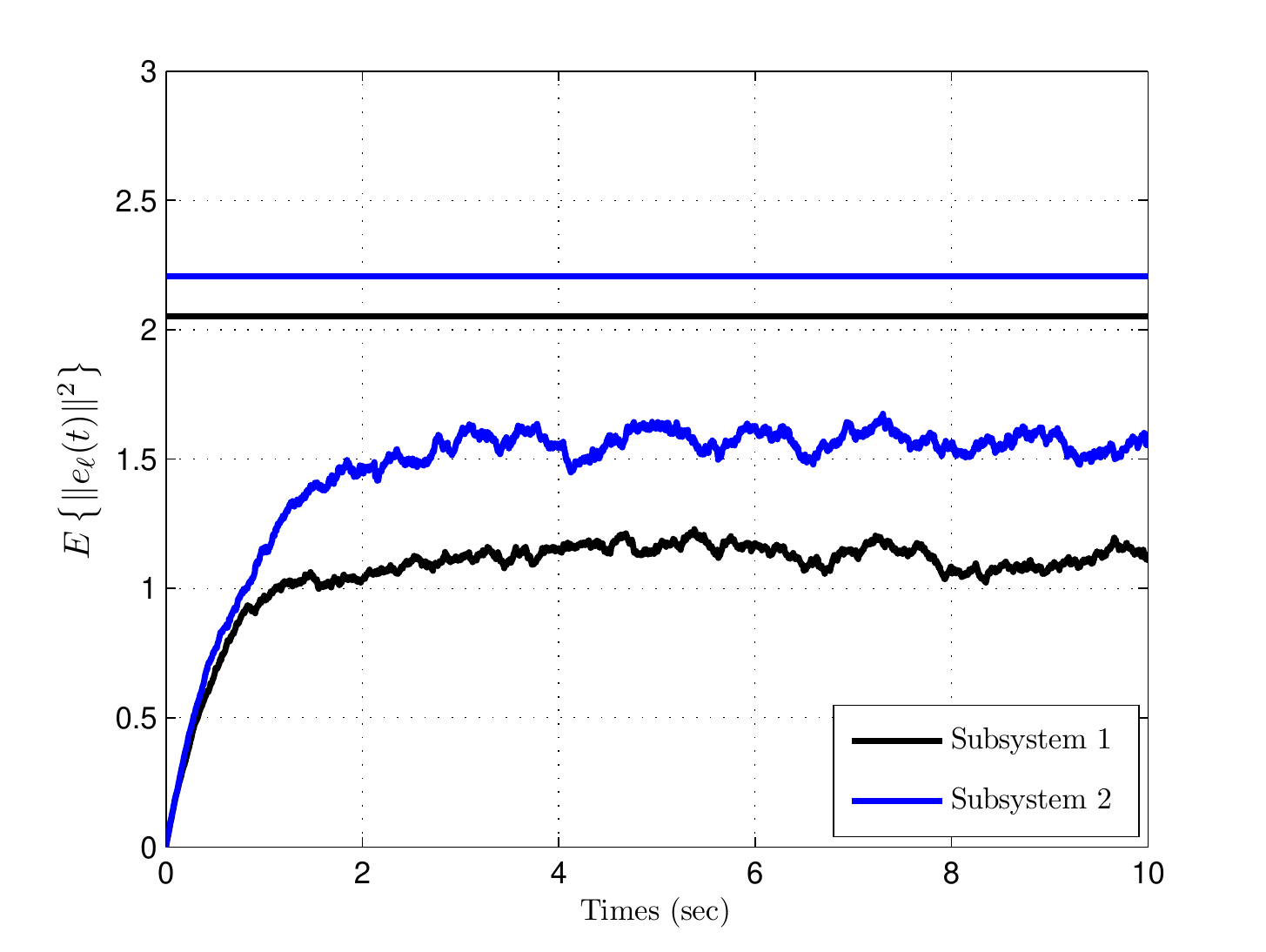}
\caption{\label{figureestimationnonscalar} Estimation error $\mathbb{E}\{\|e_\ell(t)\|^2\}$ for 1000~Monte Carlo simulations and its comparison to the theoretical results when $d_\ell=2$ for $\ell=1,2$.}
\end{figure}

\subsection{Estimation: Higher-Order Subsystems with Noisy State Measurement} \label{subsec:numerical:nonscalar:estimation}
Let us focus on a networked system composed of only two subsystems where each subsystem is a serial interconnection of two water tanks. Figure~\ref{figurewatertank}~(right) illustrates such a networked system. In this case, subsystem~$\ell$ can be described by~(\ref{eqn:generalsubsystem}) with
$$
A_\ell=\matrix{cc}{-(a_{\ell,\mathrm{T}}/a'_{\ell,\mathrm{T}})\sqrt{g/(2h_{\ell,\mathrm{T}})} & 0 \\ +(a_{\ell,\mathrm{T}}/a'_{\ell,\mathrm{T}})\sqrt{g/(2h_{\ell,\mathrm{T}})} &
-(a_{\ell,\mathrm{B}}/a'_{\ell,\mathrm{B}})\sqrt{g/(2h_{\ell,\mathrm{B}})} },
$$
where the parameters marked with $\mathrm{T}$ and $\mathrm{B}$ belong to the top and the bottom tanks, respectively. Let us fix parameters $a'_{1,\mathrm{T}}=a'_{1,\mathrm{B}}= a'_{2,\mathrm{T}} =a'_{2,\mathrm{B}}=1.00\,\mathrm{m}^2$, $a_{1,\mathrm{T}}=a_{1,\mathrm{B}}=0.20\,\mathrm{m}^2$, $a_{2,\mathrm{T}}=a_{2,\mathrm{B}}=0.10\,\mathrm{m}^2$, $h_{1,\mathrm{T}}=h_{1,\mathrm{B}}=0.40\,\mathrm{m}$, and $h_{2,\mathrm{T}}=h_{2,\mathrm{B}}=0.54\,\mathrm{m}$.

Let us assume that estimator~$\ell$ has access to the noisy state measurements of subsystem~$\ell$ (with noise variance $\mathbb{E}\{n_i^\ell n_i^{\ell \top}\}=0.09I_{2\times 2}$) at sampling instances $\{T_i^\ell\}_{i=0}^\infty$ enforced by the optimal scheduling policy described in Subsection~\ref{subsec:numerical:scalar:estimation}. Figure~\ref{figureestimationnonscalar} shows the estimation error variance $\mathbb{E}\{\|e_\ell(t)\|^2\}$. The horizontal lines in this figure show the theoretical bounds calculated in Theorem~\ref{tho:estimation:2}; i.e., $\mathbb{E}\{\|e_1(t)\|^2\} \leq 2.05$ and $\mathbb{E}\{\|e_2(t)\|^2\} \leq 2.21$. In comparison with the scalar case in Figure~\ref{figureSamplingInstance}~(lower-left), note that the bounds in Figure~\ref{figureestimationnonscalar} are less tight. The reason for this is that the dimension of the subsystems are now twice the previous case.

\begin{figure}[t]
\centering
\includegraphics[width=.4\linewidth]{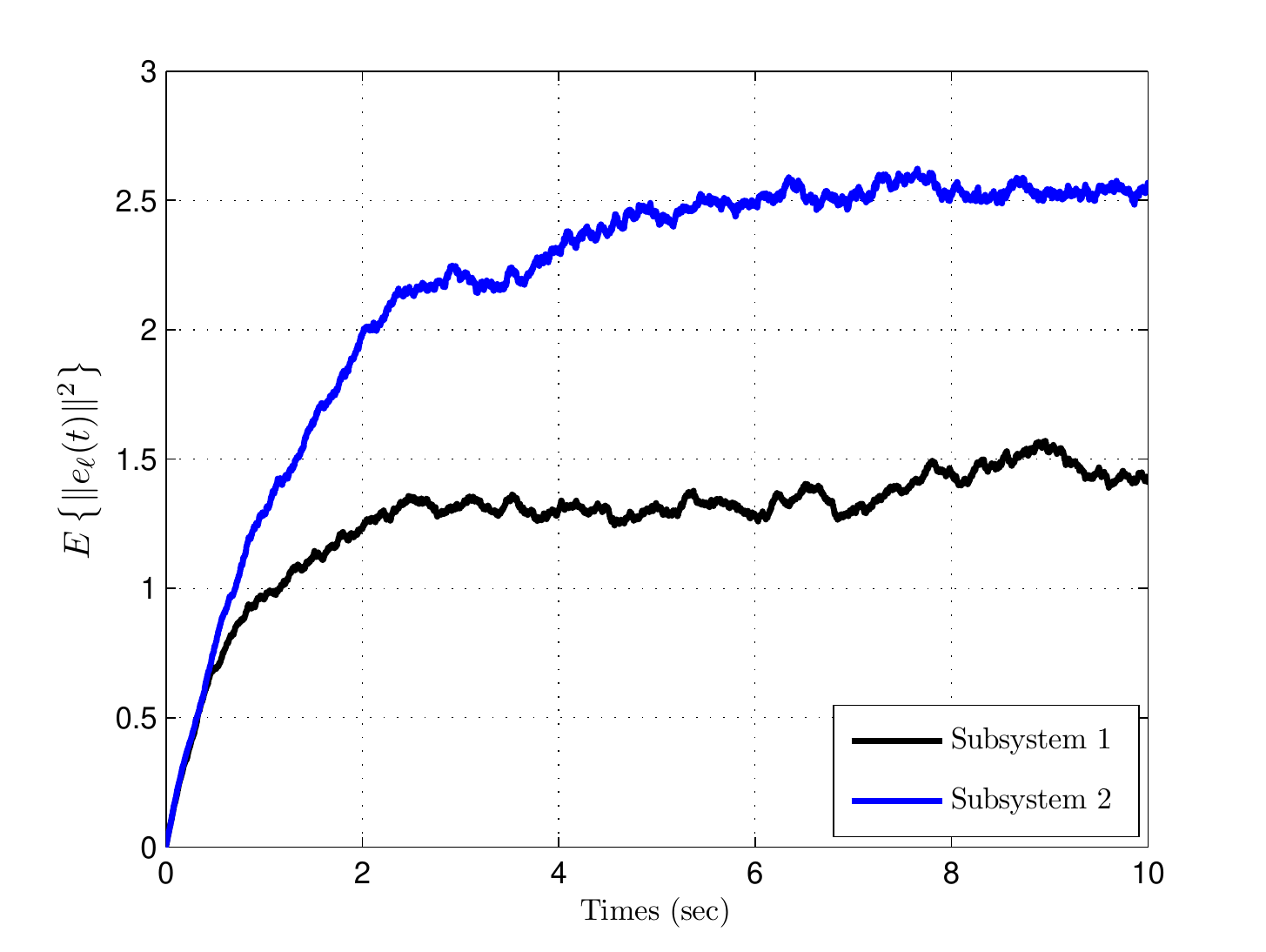}
\caption{\label{figureestimationKalman} Estimation error $\mathbb{E}\{\|e_\ell(t)\|^2\}$ for 1000~Monte Carlo simulations when using Kalman-filter based estimator. }
\end{figure}

\subsection{Estimation: Higher-Order Subsystems with Noisy Output Measurement}
In this subsection, we use output measurements $y_i^\ell=[0 \;\; 1]z_\ell(T_i^\ell)+n_i^\ell$ for all $i\in\mathbb{Z}_{\geq 0}$, where $\mathbb{E}\{(n_i^\ell)^2\}=0.09$ for $\ell=1,2$. We use the Kalman filter based scheme introduced in Subsection~\ref{subsec:Kalmanfiltering} for estimating the state of each subsystem. Figure~\ref{figureestimationKalman} illustrates the estimation error variance $\mathbb{E}\{\|e_\ell(t)\|^2\}$. As mentioned earlier, it is difficult to calculate~$\mathbb{E}\{\mathrm{trace}(P_\ell[i])\}$ as a function of the average sampling frequencies and hence, we do not have any theoretical results for comparison. Note that the upper bound presented in Theorem~\ref{tho:estimation:3} is only valid for a fixed sequence of sampling instances. This problem can be an interesting direction for future research.

\begin{figure}[t]
\centering
$$
\begin{array}{cc}
\includegraphics[width=0.4\linewidth]{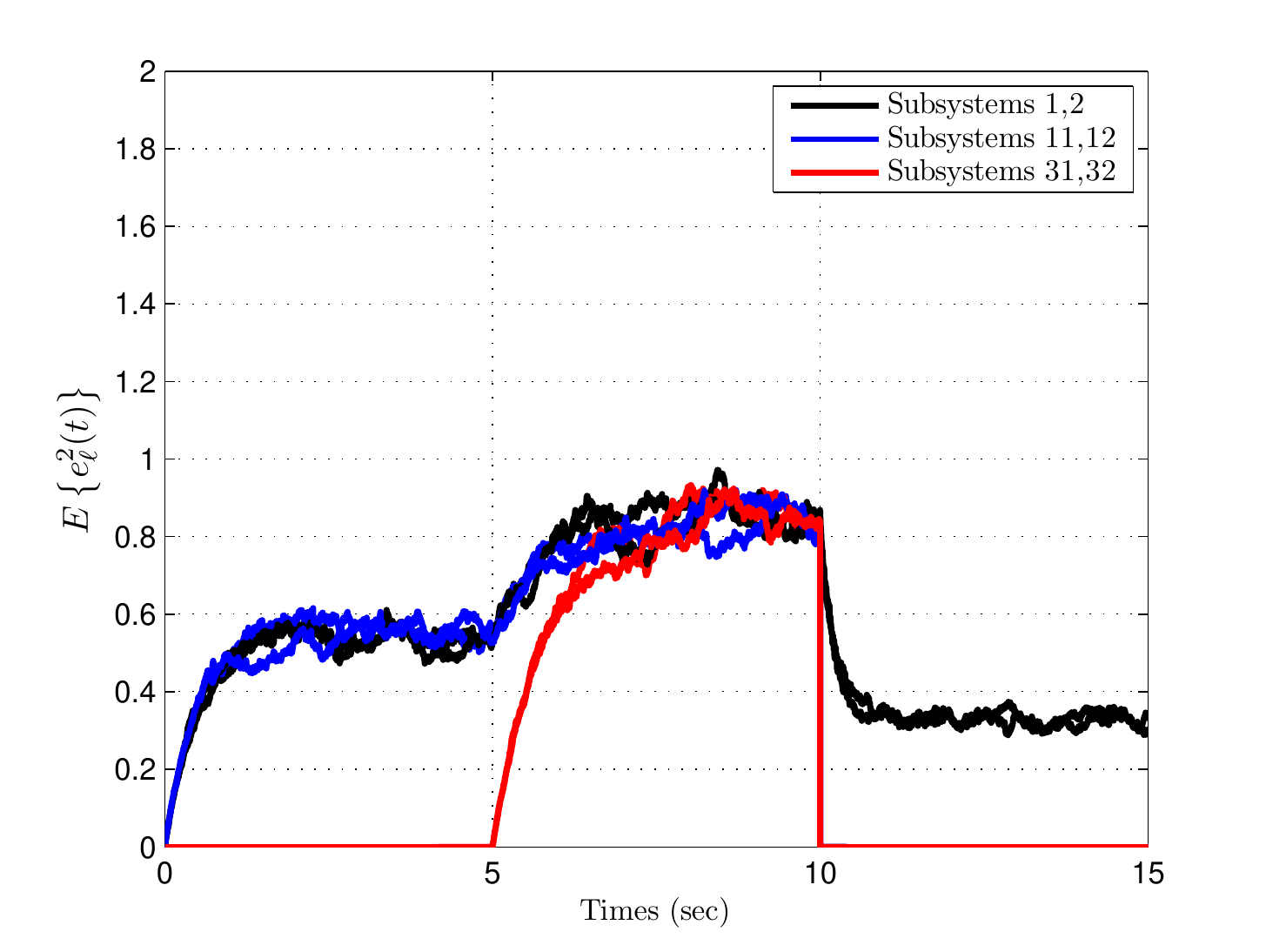} &
\vspace{-.15in}
\includegraphics[width=0.4\linewidth]{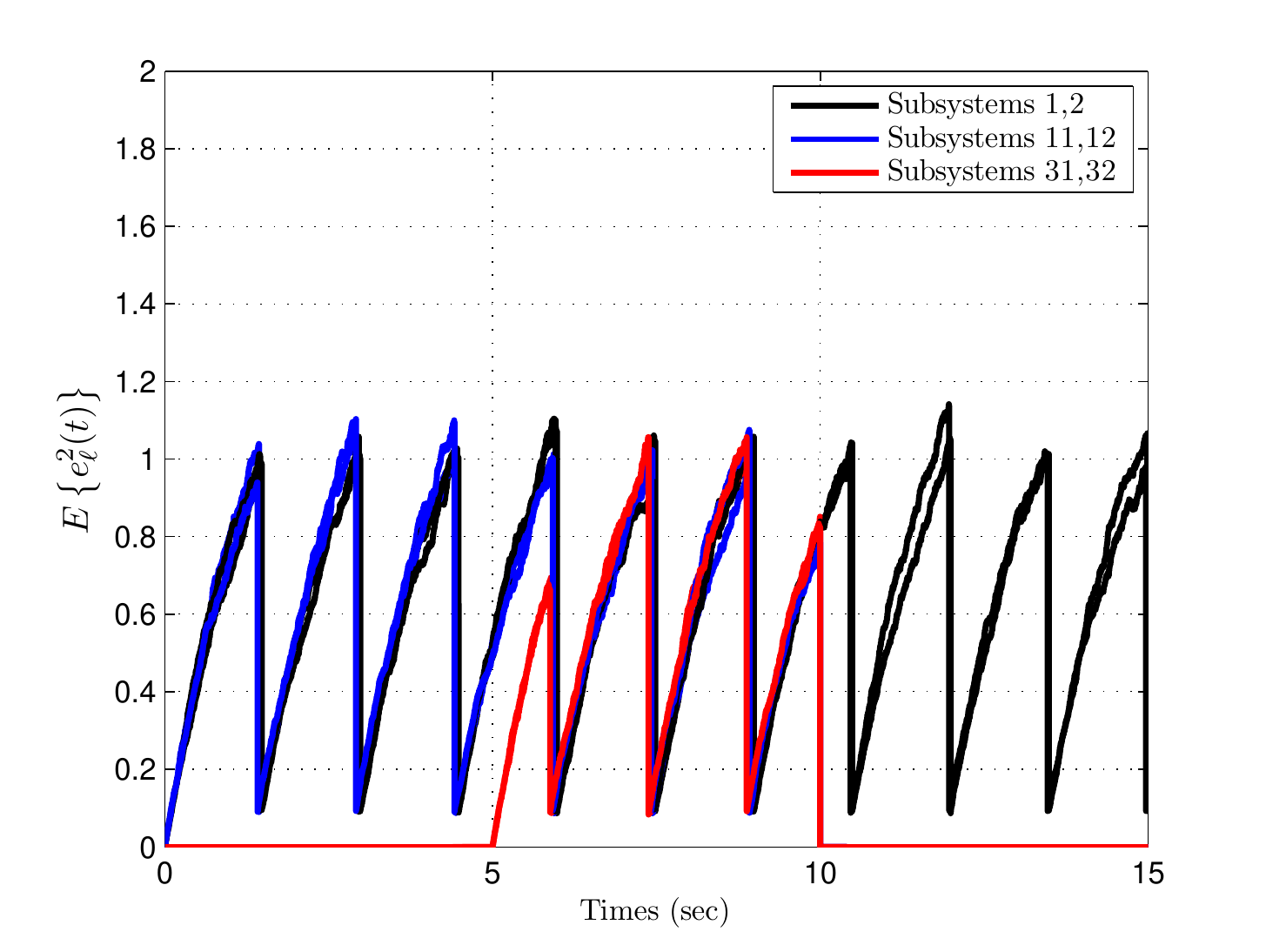}
\vspace{-.15in}
\end{array}
$$
\caption{\label{figure_adhoc_network_stochastic} Estimation error $\mathbb{E}\{e_\ell^2(t)\}$ for 1000~Monte Carlo simulations over an ad-hoc networked system with the optimal sampling policy (left) and the periodic sampling policy (right).}
\end{figure}

\subsection{Estimation: Ad-hoc Sensor Network}
As discussed earlier, an advantage of using the introduced optimal scheduling policy is that we can accommodate for changes in ad-hoc networked systems. To portray this property, let us consider a networked system that can admit up to $L=70$~identical subsystems described by~(\ref{eqn:sys}) with $\gamma_\ell=0.3$ and $\sigma_\ell=1.0$ for $1\leq \ell\leq 70$. When all the subsystems are active, we sample them using the Markov chain in~(\ref{eqn:uvr:sys}) with $m=2L=140$. We assume that $\mu_i(t)=\mu_{i,0}+u_i(t)$ for $1\leq i\leq 140$, where $\mu_{2\ell,0}=10$ and $\mu_{2\ell-1,0}=50$ for $1\leq \ell\leq 70$. In this case, we are also interested in calculating an optimal scheduling policy that minimizes
\begin{equation} \label{eqn:examplecostmanysubsystems}
J=\lim_{T\rightarrow \infty}\mathbb{E}\left\{\frac{1}{T}\int_{0}^T \sum_{\ell=1}^{70} 0.1e_{71}^\top x(t)\mathrm{d}N_{2\ell} +u(t)^\top u(t) \mathrm{d}t \right\}.
\end{equation}
However, when some of the subsystems are inactive, we simply remove their corresponding nodes from the Markov chain flow diagram in Figure~\ref{figureNCS_new}~(right) and set their corresponding terms in~(\ref{eqn:examplecostmanysubsystems}) equal to zero. Let us assume that for $t\in[0,5)$, only 30~subsystems are active, for $t\in[5,10)$, all 70~subsystems are active, and finally, for $t\in[10,15]$, only 10~subsystems are active. Figure~\ref{figure_adhoc_network_stochastic}~(left) and~(right) illustrate the estimation error variance $\mathbb{E}\{e_\ell^2(t)\}$ for 1000~Monte Carlo simulations when using the optimal scheduling policy and the periodic scheduling policy, respectively. Since in the periodic scheduling policy, we have to fix the sampling instances in advance, we must determine the sampling periods according to the worst-case scenario (i.e., when the networked system is composed of 70~subsystems).  Therefore, when using the periodic sampling, the networked system is not using its true potential for $t\in[0,5)$ and $t\in[10,15]$, but the estimation error is fluctuating substantially over the whole time range. The proposed optimal scheduling policy adapts to the demand of the system, see Figure~\ref{figure_adhoc_network_stochastic}~(right). For instance, as shown in Figure~\ref{figure_adhoc_network_stochastic}~(left), when subsystems~31 and~32 become active for $t\in[5,10)$, the overall sampling frequencies of the subsystems decreases (and, in turn, the estimation error variance increases), but  when they become inactive again for $t\in[10,15]$, the average sampling frequencies increase (and, in turn, the estimation error variance decreases). Hence, this example illustrates the dynamic benefits of our proposed stochastic scheduling approach.

\begin{figure}[t]
\centering
$$
\begin{array}{cc}
\includegraphics[width=0.4\linewidth]{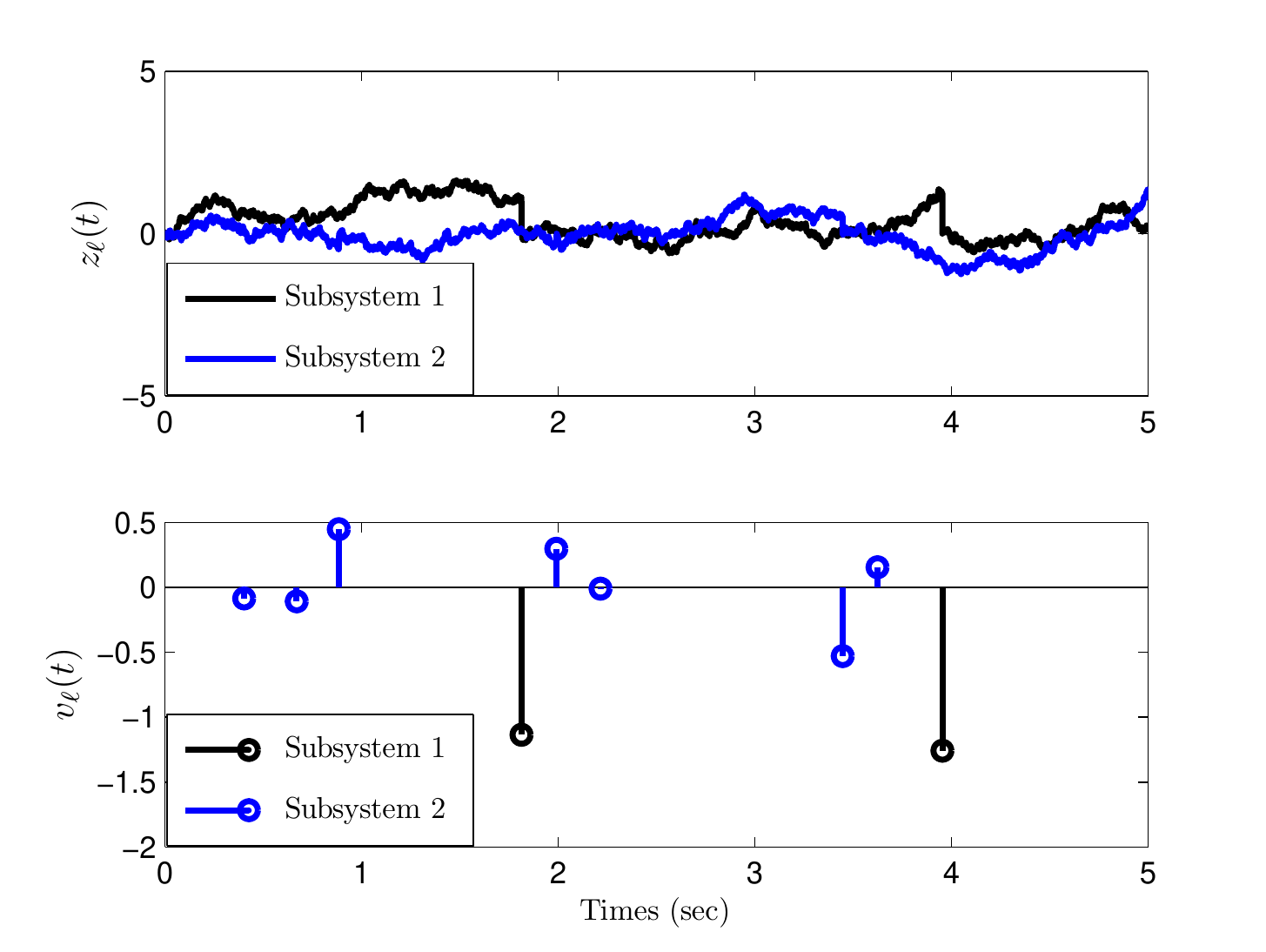} &
\vspace{-.15in}
\includegraphics[width=0.4\linewidth]{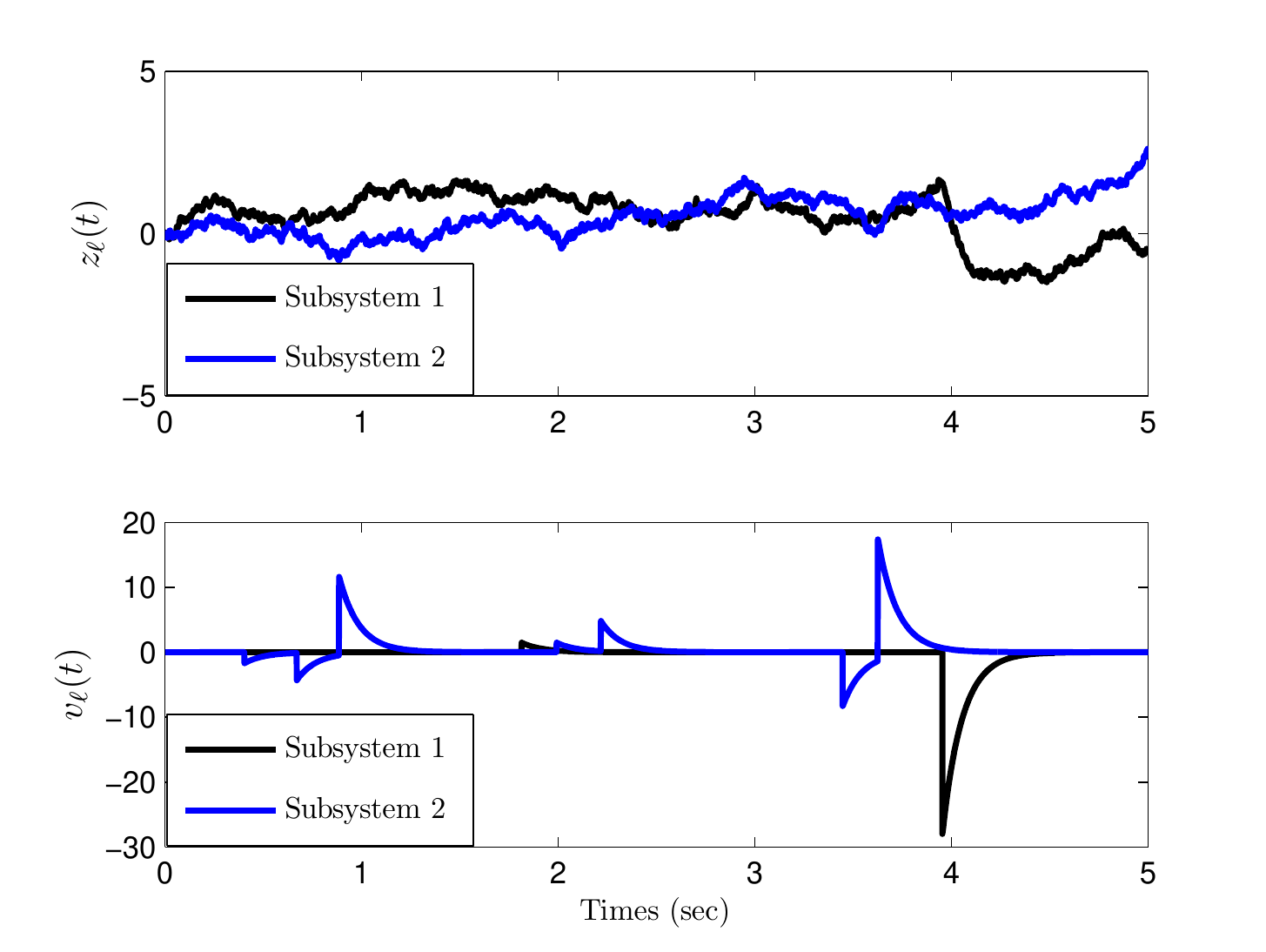}
\vspace{-.15in}
\end{array}
$$
\caption{\label{figurecontroller} An example of state and control of the closed-loop subsystems when using the impulsive controller (left) and the exponential controller with $\theta=10$ (right).}
\end{figure}

\subsection{Controller: Decoupled Scalar Subsystems}
\begin{figure}[t]
\centering
$$
\begin{array}{ccc}
\includegraphics[width=0.4\linewidth]{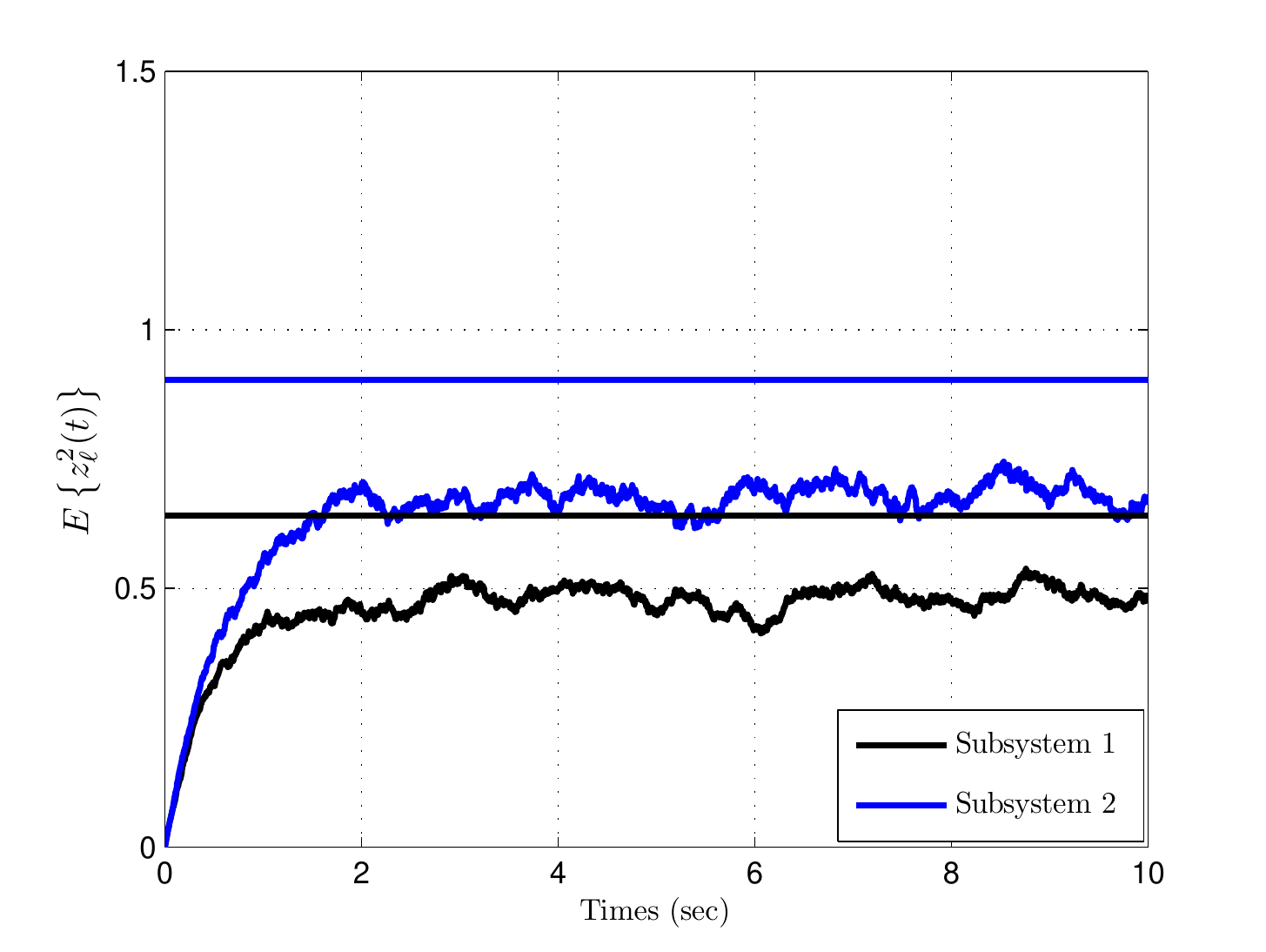} &
\vspace{-.1in}
\includegraphics[width=0.4\linewidth]{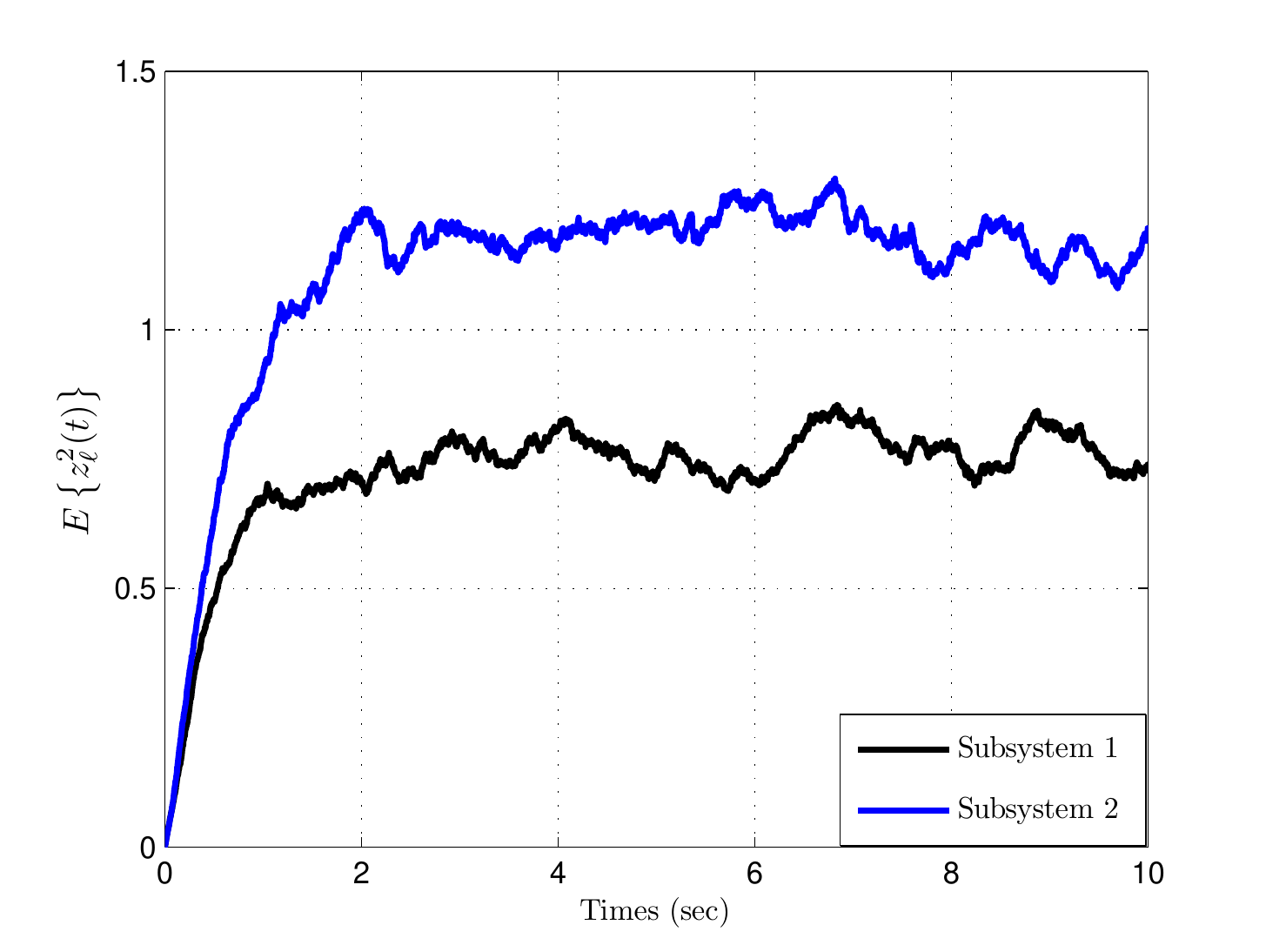} &
\vspace{-.1in}
\end{array}
$$
\caption{\label{figurecontrollerstatistics} Closed-loop performance measure $\mathbb{E}\{z_\ell^2(t)\}$ for 1000~Monte Carlo simulations when using the impulsive controller (left) and the exponential controller (right). }
\end{figure}

In this subsection, we briefly illustrate the networked control results for $L=2$ subsystems. Let the subsystems be described by~(\ref{eqn:sys:control}) with $\gamma_1=0.7$, $\gamma_2=0.3$, and $\sigma_1=\sigma_2=1.0$. Let us assume that controller~$\ell$ has access to state measurements of subsystem~$\ell$ according to~(\ref{eqn:output}) with measurement noise variance $\eta_\ell=0.3$ for $\ell=1,2$. We use the optimal scheduling policy derived in Subsection~\ref{subsec:numerical:scalar:estimation} for assigning sampling instances. Figure~\ref{figurecontroller}~(left) and~(right) illustrate an example of the state and the control signal for both subsystems when using the impulsive and exponential controllers (with $\theta=10$), respectively. Note that in Figure~\ref{figurecontroller}~(left), the control signal of the impulsive controller only portrays the energy that is injected to the subsystem (i.e., the integral of the impulse function) and not its exact value. Figure~\ref{figurecontrollerstatistics}~(left) and~(right) show the closed-loop performance $\mathbb{E}\{z_\ell^2(t)\}$ when using the impulsive and exponential controllers, respectively. The horizontal lines illustrate the theoretical upper bounds derived in Theorem~\ref{tho:control:1}. Note that the exponential controller gives a worse performance than the impulse controller. This is normal as we design the exponentials only as an approximation of the impulse train. 

\subsection{Controller: Coupled Scalar Subsystems}
\begin{figure}[t]
\centering
$$
\begin{array}{ccc}
\includegraphics[width=0.4\linewidth]{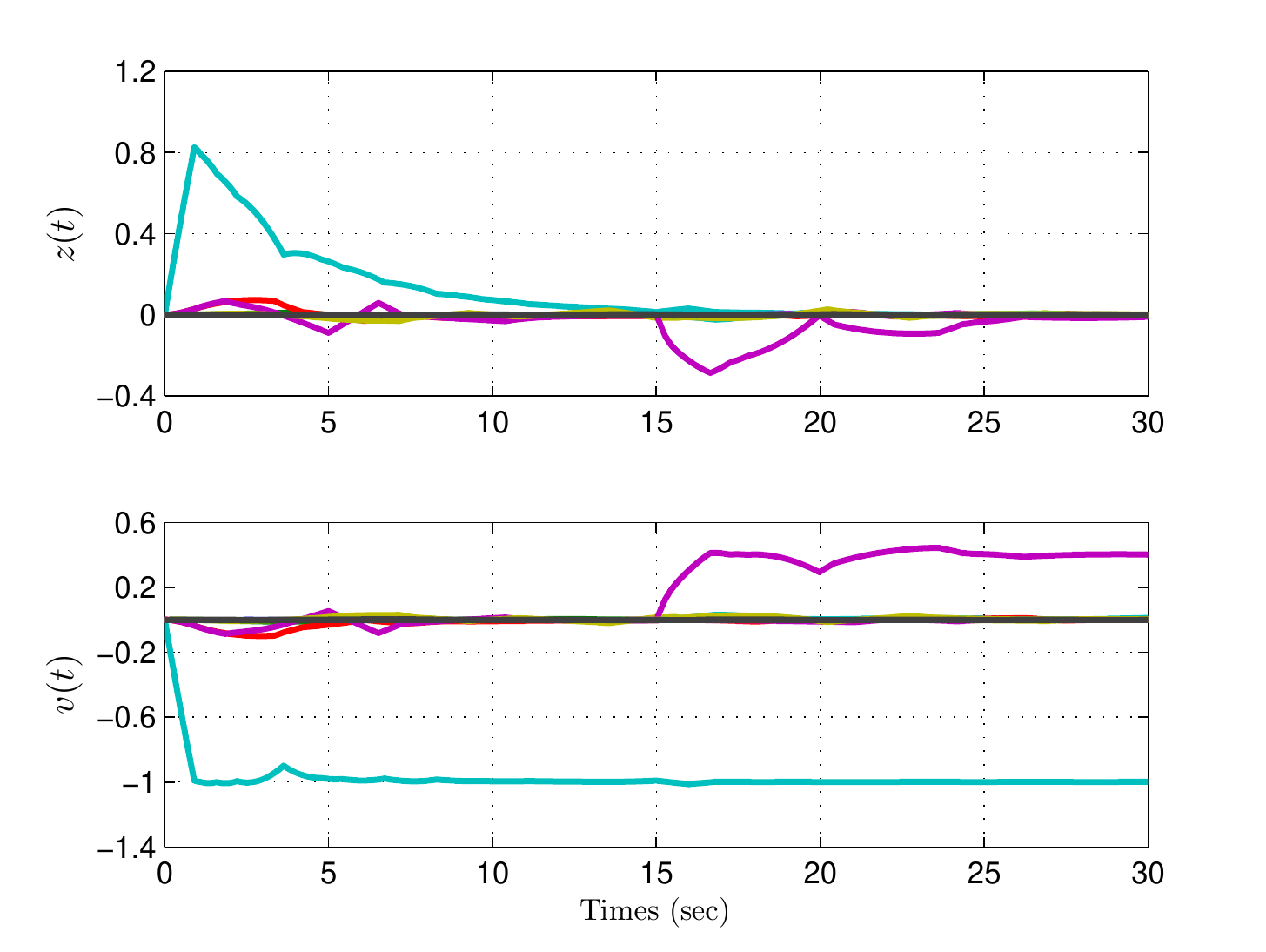} &
\vspace{-.15in}
\includegraphics[width=0.4\linewidth]{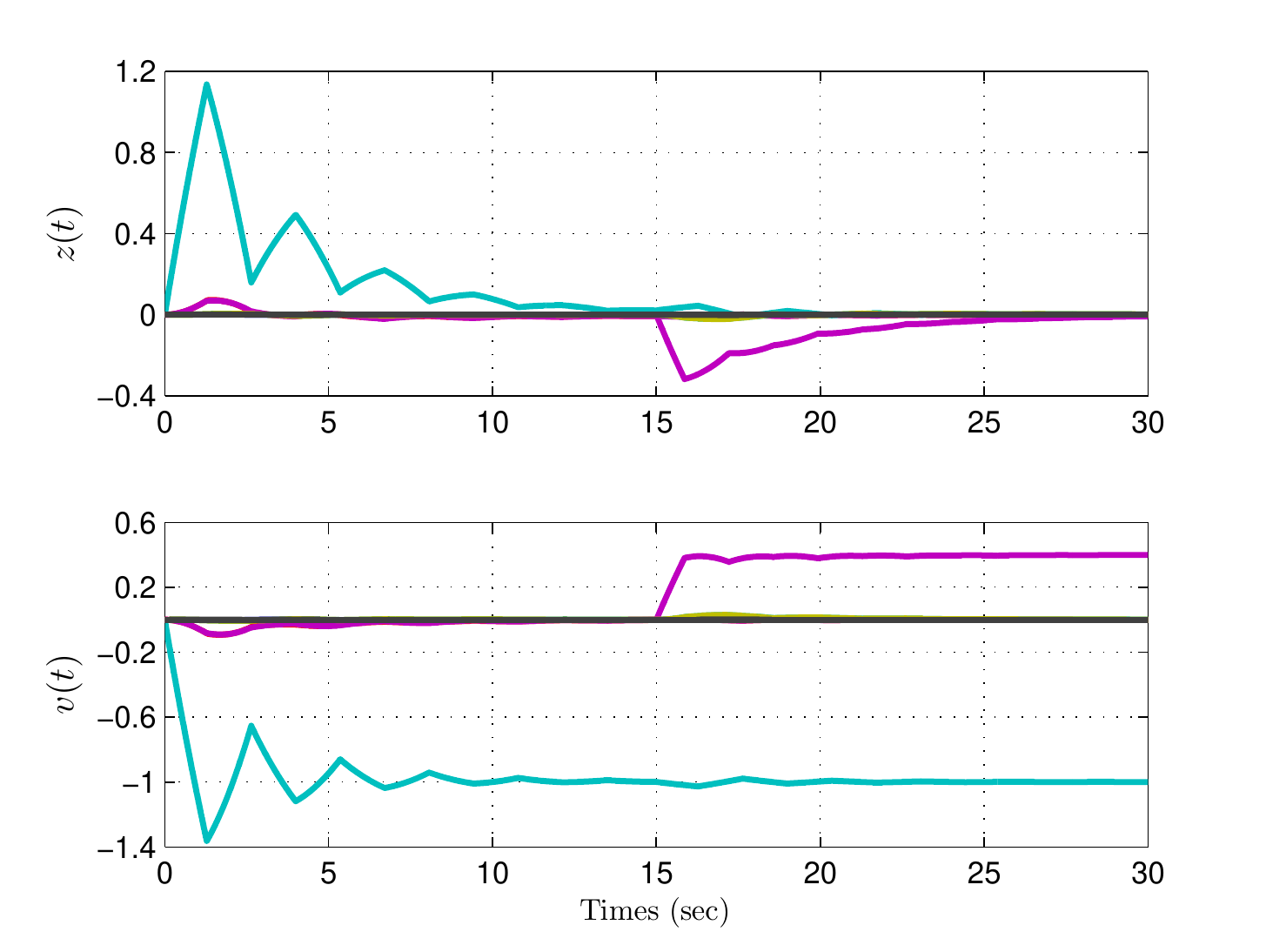}
\vspace{-.15in}
\end{array}
$$
\caption{\label{figuregeneralNCS}  An example of state and control signal using the optimal sampling policy (left) and the periodic sampling policy (right). }
\end{figure}

Consider a networked system composed of $L=70$~interconnected subsystems, where subsystem~$\ell$, $1\leq \ell\leq 70$, can be described by 
$$
\frac{\mathrm{d}}{\mathrm{d}t}z_\ell(t)=0.1(z_{\mathrm{mod}(\ell-1,L)}(t)- z_\ell(t))+0.1(z_{\mathrm{mod}(\ell+1,L)}(t)-z_\ell(t))+v_\ell(t)+w_\ell(t); \; z_\ell(0)=0.
$$
with notation $\mathrm{mod}(i,j)=i-\lfloor i/j\rfloor j$ for any $i\in\mathbb{Z}$ and $j\in\mathbb{Z}_{>0}$. In this model, $z_\ell(t)$, $v_\ell(t)$, and $w_\ell(t)$ respectively denote the state, the control input, and the exogenous input. Each subsystem transmits its state measurement over the wireless network at instances $\{T_i^\ell\}_{i=0}^\infty$ to its subcontroller. Hence, at any time $t\in\mathbb{R}_{\geq 0}$, subcontroller~$\ell$ has access to the state measurements $z_\ell(T_{M_t^\ell}^\ell)$ where recalling from the earlier definitions $M_t^\ell=\max\left\{i\geq 1\,|\, T_i^\ell\leq t \right\}$. Each subcontroller simply implement the following decentralized proportional-integral control law
$$
v_\ell(t)=-1.2z_\ell(T_{M_t^\ell}^\ell)-0.3\int_0^t z_\ell(T_{M_\tau^\ell}^\ell)\mathrm{d}\tau.
$$
We sample the subsystems using the Markov chain in~(\ref{eqn:uvr:sys}) with $m=2L=140$. We assume that $\mu_i(t)=\mu_{i,0}+u_i(t)$ for $1\leq i\leq 140$, where $\mu_{2\ell,0}=10$ and $\mu_{2\ell-1,0}=70$ for $1\leq \ell\leq 70$. Let us consider the following disturbance rejection scenario. We assume that $w_\ell(t)\equiv0$ for $\ell\neq 4,26$, $w_4(t)=\mathrm{step}(t)$, and $w_{26}(t)=-0.4\,\mathrm{step}(t-15)$, where $\mathrm{step}:\mathbb{R}\rightarrow \{0,1\}$ is the heaviside step function (i.e., $\mathrm{step}(t)=1$ for $t\in\mathbb{R}_{\geq 0}$ and $\mathrm{step}(t)=0$, otherwise). Let us denote $t\in[0,15)$ and $t\in[15,30]$ as the first phase and the second phase, respectively. During each phase, we find the infinite horizon optimal scheduling policy which minimizes
\begin{equation*}
J=\lim_{T\rightarrow \infty}\mathbb{E}\left\{\frac{1}{T}\int_{0}^T \sum_{\ell=1}^{70} \xi_\ell\, e_{71}^\top x(t)\mathrm{d}N_{2\ell} +u(t)^\top u(t) \mathrm{d}t \right\}.
\end{equation*}
We fix $\xi_\ell=10$ for $\ell=4$ over the first phase and for $\ell=26$ over the second phase. In addition, we fix $\xi_\ell=20$ for $\ell=3,5$ over the first phase and for $\ell=25,27$ over the second phase. Finally, we set $\xi_\ell=30$ for the rest of the subsystems. This way, we can ensure that we more frequently sample the subsystems that are most recently disturbed by a nonzero exogenous input signal (and the ones that are directly interacting with them). Figure~\ref{figuregeneralNCS}~(left) and~(right) illustrate an example of the system state and control input when using the described optimal scheduling policy and the periodic scheduling policy, respectively. For the periodic scheduling policy, we have fixed the sampling frequencies according to the worst-case scenario (i.e., the average frequencies of the optimal scheduling policy when $\xi_\ell=10$ for all $1\leq \ell \leq 70$ corresponding to the case where all the subsystems are disturbed). As we expect, for this particular example, the closed-loop performance is better with the optimal scheduling policy than with the periodic scheduling policy. This is indeed the case because the optimal scheduling policy adapts the sampling rates of the subsystems according to their performance requirements. 

\section{Conclusions} \label{sec:Conclusion}
In this paper, we used a continuous-time Markov chain to optimally schedule the measurement and transmission time instances in a sensor network. As applications of this optimal scheduling policy, we studied networked estimation and control of large-scale system that are composed of several decoupled scalar stochastic subsystems. We studied the statistical properties of this scheduling policy to compute bounds on the closed-loop performance of the networked system. Extensions of the estimation results to observable subsystems of arbitrary dimension were also presented. As a future work, we could focus on obtaining better performance bounds for estimation and control in networked system as well as combining the estimation and control results for achieving a reasonable closed-loop performance when dealing with observable and controllable subsystems of arbitrary dimension. An interesting extension is also to consider zero-order hold and other control function for higher-order subsystems. 

\bibliography{ref}
\bibliographystyle{ieeetr}

\end{document}